\documentclass[a4paper,12pt,oneside,reqno]{amsart}
\usepackage{hyperref}
\usepackage[headinclude,DIV=13]{typearea}
\areaset{15.1cm}{25.0cm}
\usepackage{amsfonts,amssymb,amsmath,amsthm,bbm, bm}

\usepackage[utf8]{inputenc}
\usepackage[T1]{fontenc}
\usepackage{graphicx, psfrag}
\usepackage{graphicx} 
\usepackage{color}[notcite]

\newtheorem{theorem}{Theorem}[section]
\newtheorem{lemma}[theorem]{Lemma}
\newtheorem{proposition}[theorem]{Proposition}
\newtheorem{corollary}[theorem]{Corollary}

\newtheorem{assumption}[theorem]{Assumption}

\theoremstyle{remark}
\newtheorem{remark}{Remark}

\newcommand\Bracket[2]{\genfrac{}{}{0pt}{}{#1}{#2}}

\def\paragraph#1{\noindent \textbf{#1}}

\numberwithin{equation}{section}

\def\eee{\mathrm{e}}
\def\dd{\mathtt{d}}

\def\<{\langle}
\def\>{\rangle}

\def\a{\alpha}
\def\b{\beta}
\def\e{\epsilon}

\def\g{\gamma}

\def\l{\lambda}

\def\s{\sigma}
\def\t{\tau}

\def\R{{\mathbb{R} }}  
\def\N{{\mathbb{N}}}  
\def\P{{\mathbb{P}}}

\def\C{{\mathbb{C}}} 
\def\E{{\mathbb{E}}}

\let\cal=\mathcal

\def\FF{{\cal F}}

\def\MM{{\cal M}}

\def\UU{{\cal U}}

\def\XX{{\cal X}}

\def \b {{\beta}}
\def \s {{\sigma}}

\def \t {{\tau}}

\def \g {{\gamma}}
\def \l {{\lambda}}
\def \d {{\delta}}
\def \a {{\alpha}}

\def \ba {\begin{array}}
	\def \ea {\end{array}}

\newcommand{\be}{\begin{equation}}
	\newcommand{\ee}{\end{equation}}

\newcommand{\bea}{\begin{eqnarray}}
	\newcommand{\eea}{\end{eqnarray}}
\def\TH(#1){\label{#1}}\def\thv(#1){\ref{#1}}
\def\Eq(#1){\label{#1}}\def\eqv(#1){(\ref{#1})}

\def\cov{\hbox{\rm Cov}}

\def \1{\mathbbm{1}}

\def\eee{\hbox{\rm e}}

\begin{document}
\title[The phase diagram of the CREM: Weak correlation regime.]
{
	The phase diagram of the complex
	continuous random energy model:
	The weak correlation regime.
}

\author[M.~Fels]{Maximilian Fels}
\address{M.~Fels\\Technion — Israel Institute of Technology, Faculty of Data and Decision Sciences, Haifa, 3200003,
	Israel.}
\email{felsm@campus.technion.ac.il}

\author[L.~Hartung]{Lisa Hartung}
\address{L.~Hartung\\Institut f\"ur Mathematik\\
	Johannes Gutenberg-Universit\"at Mainz\\Staudingerweg 9\\ 55099 Mainz, Germany}
\email{lhartung@uni-mainz.de}

\author[A.~Klimovsky]{Anton Klimovsky}
\address{A.~Klimovsky\\
	Institut für Mathematik\\
	Lehrstuhl für Mathematik VIII -- Angewandte Stochastik\\
	Julius-Maximilians-Universität Würzburg\\
	Emil-Fischer-Str.~30\\
	97074 Würzburg\\
	Germany
}
\email{anton.klimovsky@mathematik.uni-wuerzburg.de}



\date{\today}

\begin{abstract}
	We identify the fluctuations of the partition function of the continuous random energy model on a Galton-Watson tree in the so-called weak correlation regime. Namely, when the ``speed functions'', that describe the time-inhomogeneous variance, lie strictly below their concave hull and satisfy a certain weak regularity condition. We prove that the phase diagram coincides with the one of the random energy model. However, the fluctuations are different and depend on the slope of the  covariance function at $0$ and the final time $t$. 
\end{abstract}

\thanks{
	The work of MF was supported by a Minerva Fellowship of the Minerva Stiftung Gesellschaft fuer die Forschung mbH and partly funded by the Deutsche Forschungsgemeinschaft (DFG, German Research Foundation) 
	through Project-ID 211504053 - SFB 1060.
	The work of AK and LH was partly funded by the Deutsche Forschungsgemeinschaft (DFG, German Research
	Foundation) through Project-ID 443891315 within SPP 2265 and through Project-ID 412848929. The work of LH was partly funded by the Deutsche Forschungsgemeinschaft (DFG, German Research
	Foundation) through Projekt-ID 390685813 and GZ 2151 - Project-ID 390873048, through  Projektnummer 233630050 - TRR 146, and Project-ID 446173099. The project was initiated, when MF and LH were supported through the program ”Research in Pairs” by the Mathematisches Forschungsinstitut Oberwolfach in 2021, during the pandemic. We thank the whole team at the MFO for providing excellent working conditions and for their hospitality.
}

\maketitle

\section{Introduction}

\subsection{Variable speed BBM as CREM on a GW-tree} 

The continuous random energy model (CREM) was introduced in \cite{BK1} as a generalization of Derrida's generalized random energy model (GREM) \cite{Derrida_GREM}. Let $\mathbb{T}=(V,E)$ be a tree of depth (or time-horizon) $t$ endowed with (ultrametric) tree distance $d$. The CREM on a tree $\mathbb{T}$ is a mean zero Gaussian process $h$ on $V$ with covariances
\begin{equation}
	\label{eq:crem-correlation-structure}
	\E\left(h(v)h(w)\right)=tA(d(v,w)/t), \quad \forall v, w \in V,
\end{equation}
where
$
d(v,w)
$ denotes the time of the most recent common ancestor of $v$ and $w$
and $A\colon [0,1]\to [0,1]$ is a non-decreasing function with $A(0)=0$ and $A(1)=1.$ We choose $\mathbb{T}$ to be a continuous time supercritical Galton-Watson tree with offspring distribution $(p_k)_{k \geq 0}$ with $\sum_{k=1}^\infty p_k=1$, $\sum_{k=1}^\infty k p_k=2$ and let $K := \sum_{k=1}^\infty k(k-1)p_k<\infty$.

\begin{figure}
	\centering
	\includegraphics[width=12cm]{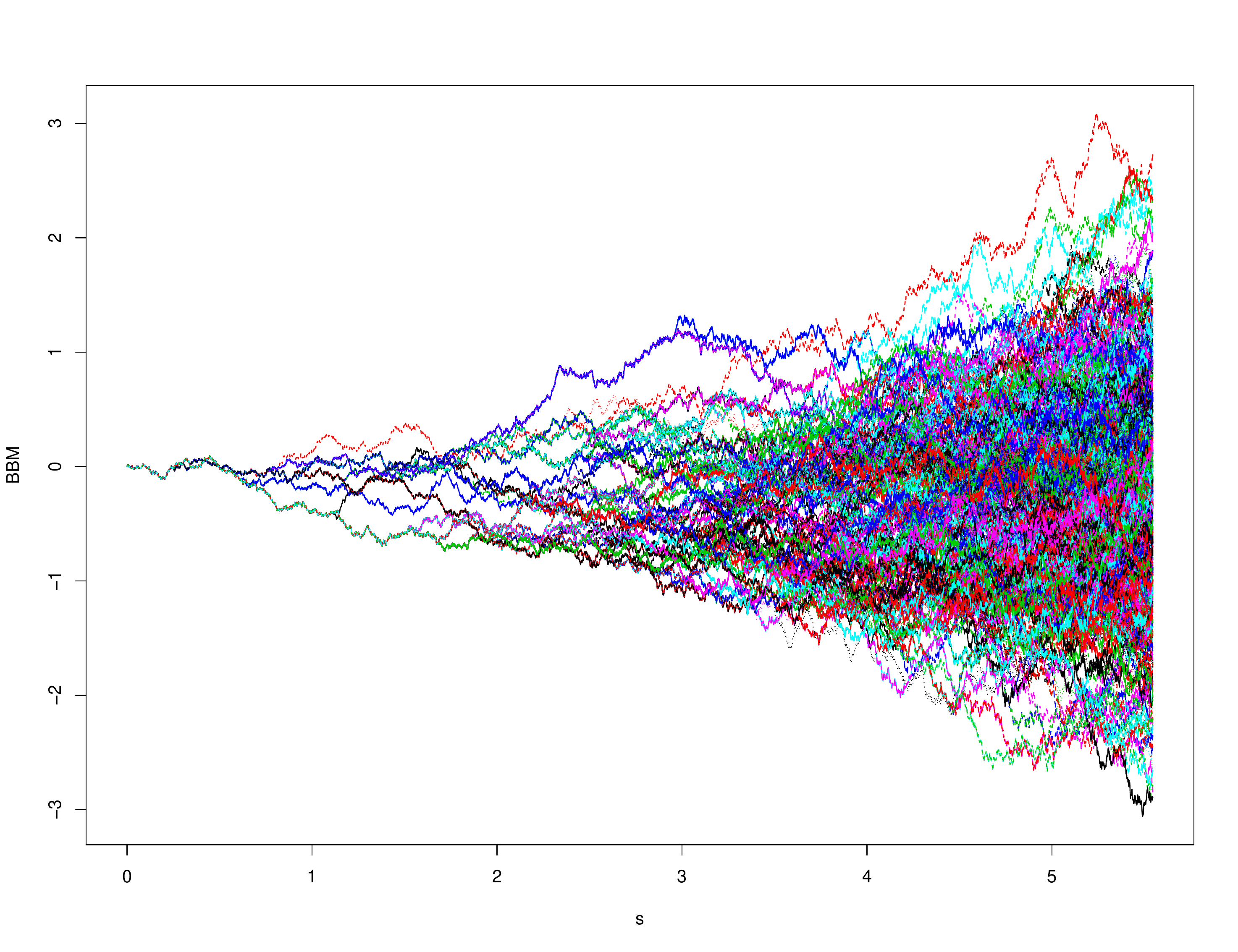}
	\caption{A simulated realization of the CREM for $A(x) := (\eee^{3x} - 1) / (\eee^{3} - 1)$ with $t = \log 256 \approx 5.55$. In this realization, $n(t) = 510$. 
		Note the increasing speed of diffusion over time.}
	\label{fig:crem}
\end{figure}

For a given time $t \geq 0$, we label the particles of the GW process as $i_1(t),\dots,
i_{n(t)}(t)$, where $n(t)$ is the total number of particles at time $t$. Note
that under the above assumptions, we have $\E\left[ n(t) \right]=\eee^t$. For $s
\leq t$ and $1\leq k \leq n(t)$, we denote by $i_k(s,t)$ the unique ancestor of particle $i_k(t)$ at
time $s$. In general, there will be several indices $k, l$ such that $i_k(s,t) =
i_l(s,t)$. 
For $t \geq 0$, the collection of all ancestors naturally induces the random tree
\begin{align}
	\label{eq:gw-tree}
	\mathbb{T}_t := \{i_k(s,t) \colon 0 \leq s \leq t, 1\leq k \leq n(t) \}
\end{align}
called the \textit{GW tree up to time $t$}. We denote by $\mathcal{F}^{\mathbb{T}_t}$
the $\s$-algebra  generated by the GW process up to time $t$. To lighten our notation, we drop the explicit dependence on $i_k(s,t)$, when referring to the particle positions. We simply denote the particle positions at time $s$ by $(x_k(s))_{k\leq n(t)}$ and mean by $x_k(s)$ the position at time $s$ of the ancestor of particle $k$ at time $t$.

\subsection{A model of complex-valued random energies}
Let  $\rho \in [-1,1]$. For any $t \in \R_+$, let $X(t) := (x_k(t))_{k\leq
	n(t)}$ and $Y(t) :=(y_k(t))_{k\leq n(t)} $ be two CREM's defined on the same underlying
GW tree such that for $k \leq n(t)$,
\begin{align}
	\label{eq:rho}
	\cov(x_k(t),y_k(t))=\rho t
	.
\end{align}
In what follows, to lighten the notation, we sometimes drop the dependence of quantities of interest on $\rho$.

We define the \textit{partition function} for the complex
BBM energy model with correlation $\rho$ at inverse temperature $\b := \s+i\t\in\C $ by
\begin{align}
	\Eq(real.1)
	\XX_{\b,\rho}(t) := \sum_{k=1}^{n(t)}\eee^{\s x_k(t)+i\t y_k(t)}
	.
\end{align}

The parameter $\rho$ specifies strength of correlation between the real and imaginary parts of the energy.
In particular, the choice $\rho = 1$ (= real and imaginary part are a.s.~the same), leads to the partition function of the real CREM at complex temperature $\beta$.



For $A(x)=\1_{1}(x)$, the model is known under the name \textit{random energy model} (REM) and the behaviour of \eqref{real.1} was analyzed in detail in \cite{KaKli14}. The case where $A$ is a step-function was analyzed in \cite{KaKli15}. Standard branching Brownian motion (corresponding to $A(x)=x$) was considered in \cite{HK15,hartklim18}. In this paper, we focus on the weak correlation regime, namely the case $A(x)<x$, for $x\in (0,1)$. 

\paragraph{Motivation and related literature.} In the physics literature, Derrida has initiated the study of random energy models (REM) as toy models of mean-field disordered systems and in particular those at complex temperatures: \cite{Derrida_zeros, DerridaEvansSpeer1993, Takahashi2011}. See \cite{MRV13, KaKli14, KaKli14_G, HK15, hartklim18} for the rigorous analysis of the REM, GREM and BBM at complex temperatures. A natural analogue in the context of log-correlated Gaussian fields is the so-called Gaussian multiplicative chaos. The case of complex temperatures has recently been treated, e.g. in \cite{LRV15, HairerShen2016, NewmanWu2019, JunnilaSaksmanWebb2020, Lacoin2023convergence}.

The motivation to consider complex partition functions is manifold, e.g., the classical
Lee-Yang theory of phase transitions \cite{LY52,bena2005statistical}, interference modelling \cite{Dobrinevski2011interference}, quantum physics \cite{peng2015experimental}, links to zeros of Riemann's $\zeta$-function (see e.g., \cite{Bailey_2022,Knauf1999-vq} for overviews). A crucial common feature in all these is the study of certain functionals of \textit{sums of complex random exponentials}. In the Lee-Yang theory, phase transitions are explained by the occurrence of the real accumulation points of the complex zeros of the so-called partition function (which is a sum of random exponenials). In quantum physics and interference, one studies the sums or integrals of complex exponentials which naturally leads to cancellations of magnitudes. Somewhat similar effects are also present in quantum spin glass models \cite{ManaiWarzel2022GREM, AdhikariBrennecke2020QuantumSK}. Concerning the Riemann $\zeta$-function, a crucial aspect of its study are striking relationships \cite{Bailey_2022}
between statistical physics of random energy models and randomized versions of the $\zeta$-function and characteristic polynomials of random matrices.

\subsection{Main results} 
In this paper, we identify the phase diagram (see  Theorem \ref{Cor:phase-diagram} below)  and  study the fluctuations of the partition function of the complex CREM in the entire subcritical regime. A natural next step is to identify the phase diagram of the complex CREM in the supercritical regime. We refer to \cite[Section~2.10]{KaKli14_G} for a conjecture on the phase diagram. Maybe surprisingly, there for strictly concave functions $A$, seven phases  are expected to emerge. The current work is a step towards understanding how to deal with the CREM at complex temperatures. 
\begin{figure}[htbp]
	\begin{minipage}[b]{0.59\textwidth}
		\centering
		\includegraphics[width=\textwidth]{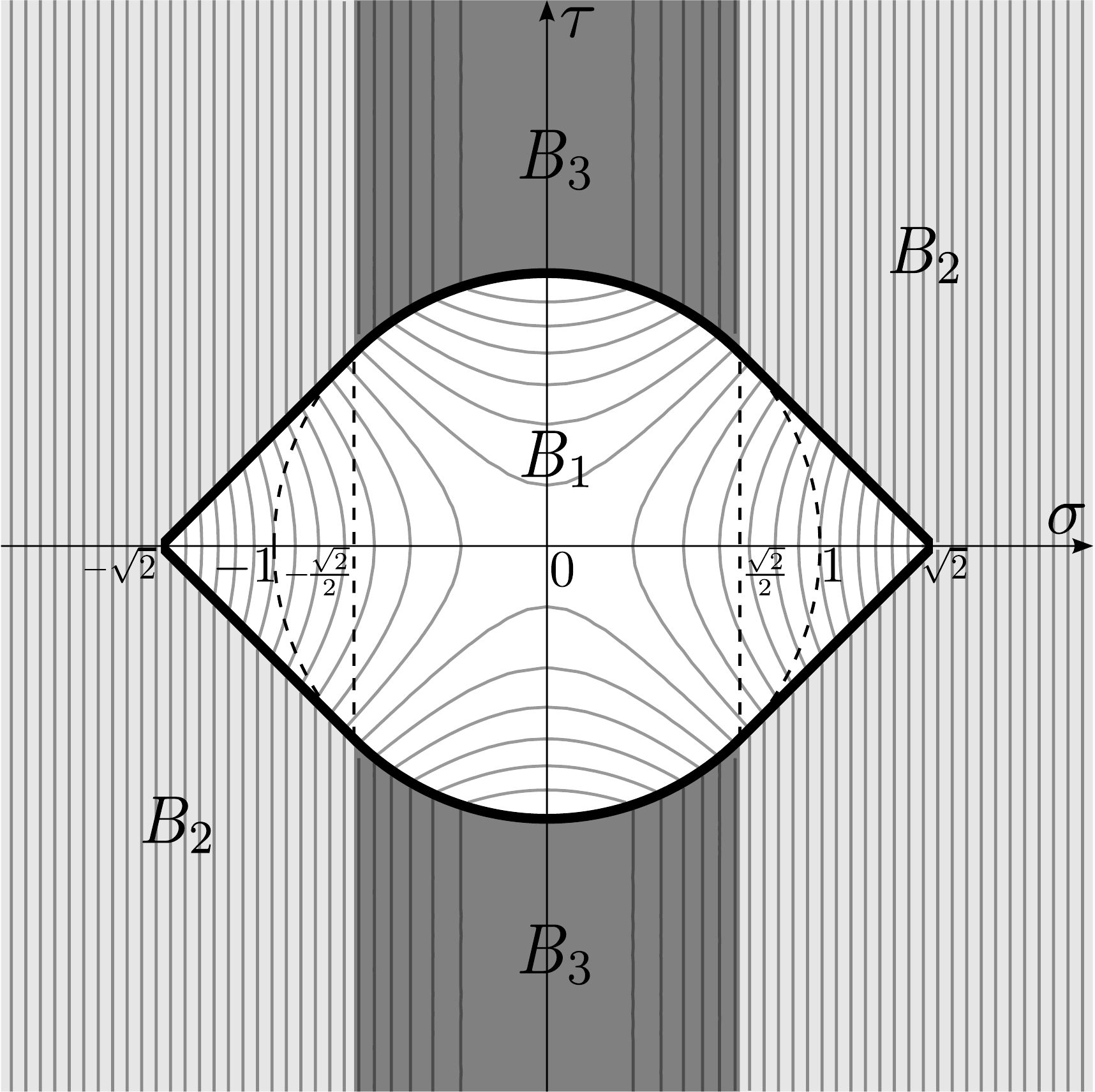}
	\end{minipage}
	\hfill
	\centering
	
	\caption{\small Phase diagram of the CREM in the weak correlation regime.
		The gray curves are the level lines of the limiting log-partition function,
		cf.~\eqref{eq:limiting-log-partition-function}. 
	}
	\vfill
	\label{fig-rem-phase-diagram}
\end{figure}

Let us specify the three domains depicted in Figure~\ref{fig-rem-phase-diagram} analytically:
\begin{equation}
	\begin{aligned}
		B_{1}
		:=
		\C \setminus \overline{B_2\cup B_3},
		\quad
		&
		B_{2}
		:=
		\{ \sigma + i \tau \in \C \colon 2\sigma^2 >  1, |\sigma|+|\tau| > \sqrt {2}\}
		,
		\\
		&
		B_{3}
		:=
		\{ \sigma + i \tau  \in \C \colon 2\sigma^2 < 1, \sigma^2+\tau^2> 1\}.
	\end{aligned}
\end{equation}

\begin{remark}
	Some of our results will be stated under the \textit{binary branching}
	assumption (i.e., $p_k = 0$ for all $k > 2$). Existence of all moments of the offspring distribution also suffices for all our results and does not require
	essential changes in the proofs.
\end{remark}
\begin{assumption} \label{weak}(Weak correlation regime, see also \cite{BH15}) Let $A\colon
	[0,1]\rightarrow [0,1]$ be a right-continuous, non-decreasing function that  satisfies the following 
	three conditions:
	\begin{itemize}
		\item[(A1)] For all $x\in (0,1)$: $A(x)<x$, $A(0)=0$ and $A(1)=1$.
		\item[(A2)]There exists $\d_b>0$ and functions $\overline{B}(x)$, $\underline{B}(x)\colon
		[0,1]\rightarrow [0,1]$ that are twice differentiable in $[0,\d_b]$ with bounded second 
		derivatives, such that
		\be\Eq(ass.2)
		\underline{B}(x)\leq A(x)\leq  \overline{B}(x),\quad \forall x\in[0,\d_b]
		\ee
		with $\overline{B}'(0)=\underline{B}'(0)\equiv A'(0)\equiv \s_b^2$.
		\item[(A3)]There exists $\d_e>0$ and functions $\overline{C}(x)$, $\underline{C}(x)\colon
		[0,1]\rightarrow [0,1]$ that are twice differentiable in $[1-\d_e,1]$ with bounded second 
		derivatives, such that
		\be\Eq(ass.3)
		\underline{C}(x)\leq A(x)\leq  \overline{C}(x),\quad \forall x\in[1-\d_e,1]
		\ee
		with $\overline{C}'(1)=\underline{C}'(1)\equiv A'(1)\equiv \s_e^2$. The case $A'(1)=+\infty$ is allowed. 
		This is to be understood in the sense that,  for all $\rho <\infty$, there exists $
		\varepsilon>0$ such that, for all 
		$x\in [1-\varepsilon,1]$, $A(x) \leq 1-\rho(1-x)$. 
	\end{itemize}

\end{assumption}
Our first result states that the complex BBM energy model indeed has the phase diagram depicted in Figure~\ref{fig-rem-phase-diagram}.

\begin{theorem}[Phase diagram]
	\label{Cor:phase-diagram}
	Let $A$ satisfy Assumption \ref{weak}. For any $\rho \in [-1,1]$ and any $\beta \in \C $, the complex CREM in the weak correlation regime and with binary branching, has the same log-partition function
	and the phase diagram (cf., \textup{Figure~\ref{fig-rem-phase-diagram}}) as the
	complex REM, i.e.,
	\begin{align}
		\label{eq:limiting-log-partition-function}
		\lim_{t \uparrow \infty} \frac{1}{t} \log \XX_{\b,\rho}(t) =
		\begin{cases}
			1 + \frac{1}{2}(\sigma^2 - \tau^2), & \beta \in \overline{B_1},
			\\
			\sqrt{2}|\sigma|, & \beta \in \overline{B_2},
			\\
			\frac{1}{2}+\sigma^2, & \beta \in \overline{B_3},
		\end{cases}
	\end{align}
	in probability.
\end{theorem}

\begin{remark} 
	(1)\,
	For a deterministic regular weighted tree (= directed polymer on a tree),
	under the assumption of no correlations between the real and imaginary parts of
	the complex random energies (i.e., case $\rho = 0$), formula
	\eqref{eq:limiting-log-partition-function} was obtained by Derrida et
	al.~\cite{DerridaEvansSpeer1993}. Our derivation of
	Theorem~\ref{Cor:phase-diagram} is based on the detailed information on the
	\textit{fluctuations} of the partition function \eqref{real.1}. 
	The arguments in \cite{DerridaEvansSpeer1993} seem to
	crucially rely on the assumption $\rho = 0$.
	\\ (2)
	It is natural to expect that the convergence in
	\eqref{eq:limiting-log-partition-function} also holds in $L^1$, see
	\cite[Theorem~2.15]{KaKli14} for a related result for the REM.
\end{remark}

\paragraph{Phase $\boldsymbol B_{\boldsymbol 1}$.}
We define the following martingales
\begin{equation}\Eq(mart.1)
	\MM_{\s,\t}(t) := 
	\sum_{k=1}^{n(t)} \eee^{-t\left(1+2i\rho\s\t+\frac{\s^2-\t^2}{2}\right)}\eee^{\s \tilde{x}_k(t)+i\t \tilde{y}_k(t)},
\end{equation}
where $\{\tilde{x}_k(t), k\leq n(t)\}, \{\tilde{y}_k(t), k\leq n(t)\}$ are standard BBMs ($A(x)=x$) with correlation constant $\rho$.

It was shown in \cite{hartklim18} that the limit
\be\Eq(B1.1) \lim_{t\uparrow
	\infty}\mathcal{M}_{\s,\t}(t) =: \mathcal{M}_{\s,\t} 
\ee 
exists a.s., in
$L^1$, and is non-degenerate.


\begin{remark}\label{rem:coupling}
	One can couple a realization of $(x_k(r))_{k\leq n(r)}$ to a realization of $(\sigma_b \tilde x_k(r))_{k\leq n(r)}$, both defined on the same underlying supercritical Galton-Watson tree:
	To each edge in the tree we associate an independent Gaussian with mean zero and variance equal to the length of the edge.
	We get a version of $(\sigma_b \tilde x_k(r))_{k\leq n(r)}$ by multiplying each Gaussian by $\sigma_b$.
	To obtain $(x_k(r))_{k\leq n(r)}$ using the same Gaussian random variables, we multiply each Gaussian by the square root of the corresponding increment in the covariance function. 
\end{remark}

\begin{theorem}\label{paseB1}  Let $A$ satisfy Assumption \ref{weak}.
	Let $\beta =\s+i\t$ with $\beta\in B_1$, and $\rho\in[-1,1]$ and let $\widehat{ \mathcal{M}}_{\s_b\s,\s_b\t}$ be the coupled realization of $\mathcal{M}_{\s_b\s,\s_b\t}$. Then,
	\be
	\eee^{-t\left(1+2i\rho\s\t+\frac{\s^2-\t^2}{2}\right)}\XX_{\b,\rho}(t) \to \widehat{\mathcal{M}}_{\s_b\s,\s_b\t},\qquad \text{as }t\uparrow \infty,
	\ee
	in probability.
\end{theorem}
\paragraph{Phase $\boldsymbol B_{\boldsymbol 2}$.} 
In phase $B_2$, the behaviour of the partition function is influenced by the extremes of $X$. In the weak correlation regime the extremal process is well understood
\cite{BovHar13, BH15}. Namely, for $m(t)=\sqrt{2} t-\frac{3}{2\sqrt{2}}\log t$ and as $t$ tends to infinity, 
\be\label{eq:extremal_process}
\sum_{k=1}^{n(t)}\d_{x_k(t)-m(t)}\Rightarrow 
\sum_{k\in \N}\sum_{l \in \N}\d_{\eta_k+\Delta^{(k)}_l}\quad \text{in law},
\ee
where $\eta_k$ are the atoms of a Cox process with intensity $C(\s_e)Y_{\s_b} e^{-\sqrt{2}  x}$, where $Y_{\s_b}$ is a random variable depending only on $\s_b=A'(0)$ and $C(\s_e)$ is a constant depending only on $\s_e=A'(1)$. $(\Delta^{(k)}_l)_{l\in \N}$ are the atoms of independent copies of the point process $\s_e\Delta$, where $\Delta$ is defined as the limit (in law) of
\be 
\sum_{k=1}^{n(t)}\d_{\bar x_k(t)- \max_{k\leq n(t)}\bar x_k(t) }, \mbox{ conditioned on } \max_{k\leq n(t)}\bar x_k(t)\geq \sqrt{2}\s_e t,
\ee
as $t\to \infty$. Here, $(\bar x_k(t))_{k\leq n(t)}$ denotes the particle positions of a standard BBM.

In the phase $B_2$, the limiting partition function can be described as follows.

\begin{theorem}\label{phaseB2.1} Let $A$ satisfy Assumption \ref{weak}.
	For $\beta=\s+i \t$ with $\beta \in B_2$ and $|\rho|=1 $, the rescaled partition function $e^{-\beta m(t)}\XX_{\b,\rho}(t)$ converges in law to the random variable
	\be
	\XX_{\b,1}:=\sum_{k,l\geq 1}e^{\beta \eta_k+\Delta^{(k)}_l},\qquad \text{as }t\uparrow \infty,
	\ee
	where $\eta_k$ and $\Delta^{(k)}_l$ are as in \eqref{eq:extremal_process}.
\end{theorem}
\begin{theorem}\label{phaseB2.2} Let $A$ satisfy Assumption \ref{weak}.
	For $\beta=\s+i \t$ with $\beta \in B_2$ and $|\rho|\in (-1,1) $, the rescaled partition function $e^{-\s m(t)} \XX_{\b,\rho}(t)$ converges in law as $t\uparrow \infty$ to the random variable 
	\begin{align}
		\XX_{\b,\rho}=\sum_{k,l\geq 1}e^{\sigma(\eta_k+\Delta^{(k)}_l)}U^{(k)}\tilde{W}^{(k)}_l,
	\end{align} 
	where $(U^{(k)})_{k\geq 1}$ are i.i.d.~uniformly distributed on the unit circle, $(\tilde{W}^{(k)}_l)_{l\geq 1}$, for $k\geq 1$, the atoms of independent point processes on the unit circle and where $\eta_k$ and $\Delta^{(k)}_l$ are as in \eqref{eq:extremal_process}. Moreover, the law of $\XX_{\beta,\rho}$ is complex isotropic $\sqrt{2}/\s-$stable.
\end{theorem}
\paragraph{Phase $\boldsymbol B_{\boldsymbol 3}$.} To state the next result, we need additional notation. We denote by
$\mathcal{L}[\cdot]$, $\mathcal{L}[\cdot \mid \cdot ]$, and $\Longrightarrow$, the law, conditional law, and weak convergence
respectively. By $\mathcal{N}(0, s^2)$, $s^2 > 0$, we denote the centered complex isotropic Gaussian distribution with density 
\begin{align}
	\C \ni z \mapsto \frac{\eee^{-|z/s|^2}}{\pi s^2} \in \R_+
\end{align}
w.r.t.~the Lebesgue measure on $\C$.

\begin{theorem}[CLT with random variance in $B_3$]\TH(clt.B3)Let $A$ satisfy Assumption~\ref{weak}.
	For $\b \in B_3$, $\rho\in[-1,1]$ and binary branching, e.g., $p_k=1$ for $k=2$ and $p_k=0$ otherwise. Then, as $t\uparrow \infty$,
	\begin{equation}\Eq(clt.1.1)
		\mathcal{L}
		\left[
		\frac{\XX_{\b,\rho}(t)}{ \eee^{t(1/2+\s^2)}} ~\Big\vert~ \MM_{2\s\s_b,0}
		\right]
		{\Rightarrow} \mathcal{N}\left(0,C_2\MM_{2\s\s_b,0}\right)
		,
	\end{equation}
	where $C_2 > 0$ is some constant and convergence is in law.
	
\end{theorem}
\begin{remark}
	A CLT similar to the one in \eqref{clt.1.1} should hold in $B_1$ as long as $|\s|<1/\sqrt{2}$. The proof for standard BBM in \cite{hartklim18} checks a Lindeberg-Feller condition. Alternatively, one could use the methods of moments as is done in the proof of Theorem \thv(clt.B3). We chose to not include this into the present article as it would not give many more insights. 
\end{remark}
\paragraph{Outline of the article.} In Section 2, we prove a general upper envelope for the paths of particles in the CREM. In Section 3, we study the partition function in Phase $B_1$ and prove Theorem \ref{paseB1}. Section 4 deals with Phase $B_2$ and Section 5 with Phase $B_3$. Finally, we prove Theorem \ref{Cor:phase-diagram} in Section 6.



\section{Upper envelope}
In this subsection, we provide a sufficiently tight upper envelope for the particles up to a fixed time $t>0$.
Let
\be\Eq(B1.111)
m_A(s)= \sqrt{2 s t A(s/t)} -\frac{\sqrt{tA(s/t)}}{2\sqrt{2s}}\log(s).
\ee
Note the implicit dependence of $m_A$ on $t$.
Since $A(1)=1$, we get for $s = t$ in \eqref{B1.111} that
\begin{equation}
	m_A(t)\equiv m(t)= \sqrt{2} t -\frac{1}{2\sqrt{2}} \log(t).
\end{equation} 
Fix a parameter $\g > 0$ and let $C>0$ be a constant (chosen later to be sufficiently large). We introduce the following deterministic \textit{upper envelope}
\be\Eq(B1.3)
U_{A,\g}(s) := m_A(s) + \big(((tA(s/t) \wedge (t-A(s/t)t))\vee C\big)^\g, \quad s \in [0,t],
\ee
and consider the following subset of paths of a particle respecting the upper envelope:
\be\Eq(B1.2)
\UU_{t,\g} : =\left\{x(\cdot)\in C\left([0,t],\R\right) : x(s)\leq  U_{A,\g}(s),  \forall s \in [0,t]\right\}
.
\ee
The next lemma shows that overshooting the upper envelope on $[0,t]$ is unlikely for the system of branching particles $\{x_k \colon [0,t] \to \R \}_{k \leq n(t)}$. 
\begin{lemma}\label{Lem.B1.1}
	Let $\g<1/2$. Then, for any $\e>0$, there exists $C>0$ such that for all $t>0$ sufficiently large,
	\be\Eq(B1.4)
	\P\left(\exists k\leq n(t) \colon \exists s \in \left[0,t\right] : x_k(s)>U_{A,\g}(s)\right)<\e.
	\ee
\end{lemma}
We adapt the proof in the case of usual BBM (cf., \cite[Theorem~2.2]{ABK_G}).
We first prove that it is very unlikely for the maximum of the process to cross the upper envelope $U_{A,\gamma}$ at integer times. In a second step, we extend this to all times using Gaussian estimates.

\begin{lemma}\label{Lem.B1.1int}
	Let $\g<1/2$. Then, for any $\e>0$, there exists $C>0$ such that for all $t$ sufficiently large,
	\begin{align}\label{B1.4.1}
		\P\left(\exists k\leq n(t) \mbox{ s.t.\ } \exists s\in \left[0,t\right] : x_k(\lceil s \rceil)>U_{A,\g}(\lceil s \rceil )\right)<\e.
	\end{align}
\end{lemma}
\begin{proof}
	Let
	\begin{equation}
		r=r(C)=\min\left\{\Bracket{\sup \{s\in [0,t]: U_{A,\g}(t-s)= m_A(t-s)+C^\g \},}{\sup \{s\in [0,t]: U_{A,\g}(s)= m_A(s)+C^\g \}}\right\}.
	\end{equation} The event
	\begin{equation}
		\left\{\exists k\leq n(t) : \exists s \in \left[0,t\right] : x_k(\lceil s \rceil)>U_{A,\g}(\lceil s \rceil ) \right\}
	\end{equation} 
	is the union over $j \in [t]$ of the events
	\begin{align}\label{B.1.4.2}
		\left\{\exists k\leq n(t): x_k(j)>U_{A,\g}(j) \right\}.
	\end{align}
	By a union bound and Markov's inequality, we get
	\begin{align}\label{eq:2.9}
		&\P\left(\bigcup_{j=1}^{\lfloor t\rfloor}\exists k\leq n(t): x_k(j)>U_{A,\g}(j)  \right)
		\leq \sum_{j=1}^{\lfloor t \rfloor} \eee^j \P\left(x(j)>U_{A,\g}(j)  \right),
	\end{align}
	where $x(j)$ is a Gaussian random variable with mean zero and variance $tA(j/t)$. Here, we may assume without loss of generality that $A(x)>0$ for $x>0$, otherwise starting with $j=\lfloor\sup\{s\in [0,t]: A(x)=0\}\rfloor+1$. We bound \eqref{eq:2.9} using Gaussian tail estimates. 
	We distinguish the cases $tA(j/t)\leq t-A(j/t)t$ and $tA(j/t)>t- A(j/t)t$. Let $j_1=\arg\max_{1\leq j \leq t }{tA(j/t)\leq t-A(j/t)t}$. \eqref{eq:2.9} is bounded from above by
	\begin{align}
		&\sum_{j=1}^{\lfloor (\Sigma^2)^{-1}(r) \rfloor }\frac{\eee^j \sqrt{tA(j/t)}}{(m_A(j)+C^\g)\sqrt{2\pi}} \exp\left[-\frac{\left(m_A(j)+C^\gamma\right)^2}{2tA(j/t)}\right]\nonumber\\
		& \quad+\sum_{j=\lfloor (\Sigma^2)^{-1}(r) \rfloor}^{\lfloor j_1 \rfloor}
		\frac{\eee^j \sqrt{tA(j/t)}}{(m_A(j)+(tA(j/t))^{\g})\sqrt{2 \pi}}\exp\left[-\frac{\left(m_A(j)+(tA(j/t))^{\g} \right)^2}{2tA(j/t)} \right]\nonumber
		\\&\quad+  \sum_{j=\lfloor j_1 \rfloor+1}^{\lceil (\Sigma^2)^{-1}(t-r)\rceil}
		\frac{\eee^j\sqrt{tA(j/t)}}{(m_A(j)+(t-tA(j/t))^{\g})\sqrt{2 \pi}}\exp\left[-\frac{\left(m_A(j)+(t-tA(j/t))^{\g} \right)^2}{2tA(j/t)} \right]\nonumber\\
		&\quad +\sum_{j=\lceil (\Sigma^2)^{-1}(t-r)\rceil+1}^{t}\frac{\eee^j \sqrt{tA(j/t)}}{(m_A(j)+C^\g)\sqrt{2\pi}} \exp\left[-\frac{\left(m_A(j)+C^\gamma\right)^2}{2tA(j/t)}\right].
	\end{align}
	We observe that since $A(x)<x$, $r(C)\leq C$. Moreover, since $m_A(j)$ is the order of the maximum of $\eee^j$ i.i.d.~centered Gaussian random variables with variance $tA(j/t)$, one can easily verify that the right-hand side of \eqref{eq:2.9} is bounded from above by a constant times
	\begin{align}\label{eq:2.10}
		2C\exp\left[-\sqrt{2}C^\g] \right]&+
		\sum_{j=\lfloor (\Sigma^2)^{-1}(r) \rfloor}^{\lfloor j_1 \rfloor}
		\exp\left[-\sqrt{2} (tA(j/t))^{\g}\right]\nonumber\\
		&+  \sum_{j=\lfloor j_1 \rfloor+1}^{\lceil (\Sigma^2)^{-1}(t-r)\rceil}
		\exp\left[-\sqrt{2}\left(t-tA(j/t) \right)^{\g} \right].
	\end{align}
	The last two sums in \eqref{eq:2.10} are finite and tend to zero as $t\rightarrow \infty$ and so the probability in \eqref{B1.4.1} can be made arbitrarily small when taking limits $t\rightarrow \infty$ followed by $C\rightarrow \infty$, which implies the claim.
\end{proof}
Lemma~\ref{Lem.B1.1int} allows deducing Lemma~\ref{Lem.B1.1}, which extends it to continuous time. This is done using the fact that, if the event of exceeding the upper envelope does happen at time $s$, the maximum at time $\lceil s \rceil$ is very likely to be still high, namely greater than $U_{A,\gamma/2}(\lceil s \rceil)$. By Lemma~\ref{Lem.B1.1int}, the probability of the latter event can be made arbitrarily small by choosing $C>0$ sufficiently large.
\begin{proof}[Proof of Lemma~\ref{Lem.B1.1}:]
	Considering the cases whether the maximum of BBM at time $j$ is smaller, equal or larger than $U_{A,\g}(s)$, we bound the probability of the event,
	\begin{align}
		\left\{\exists k\leq n(t) \colon \exists s \in \left[0,t\right] : x_k(s)>U_{A,\g}(s)\right\}
	\end{align}
	from above by
	\begin{align}\label{eq:2.12}
		\P&\left(\exists k\leq n(t) \colon \exists s \in \left[0,t\right] : x_k(s)>U_{A,\g}(s)\land x_k(\lceil s \rceil )>U_{A,\g / 2}(\lceil s \rceil ) \right) \nonumber\\
		&+ \P\left(\exists k\leq n(t) \colon \exists s \in \left[0,t\right] : x_k(s)>U_{A,\g}(s)\land x_k(\lceil s \rceil)<U_{A,\g / 2}(\lceil s \rceil ) \right). 
	\end{align}
	The first probability in \eqref{eq:2.12} is bounded from above by
	\begin{align}
		\P\left( \exists k\leq n(t) \colon \exists s \in \left[0,t\right] : x_k(\lceil s \rceil )>U_{A,\g / 2}(\lceil s \rceil ) \right).
	\end{align}
	By Lemma~\ref{Lem.B1.1int}, this is bounded from above by $\epsilon/2$ for $C>0$ large. It remains to bound the second probability. Let $\mathcal{T}$ denote the stopping time
	\begin{align}
		\mathcal{T}\equiv \inf\left\{ s\in \left[0,t\right]:\, \exists k\leq n(t) : x_k(s)>U_{A,\g}(s) \right\}.
	\end{align}
	By conditioning on $\mathcal{T}$, we can rewrite the second probability in \eqref{eq:2.12} as
	\begin{align}\label{eq:2.15}
		\int_{0}^{t}\P\left(\exists k\leq n(t) : x_k(s^\prime)<U_{A,\g/2}(s^\prime) \mid \mathcal{T}=s^\prime\right) \P\left( \mathcal{T}\in \mathrm{d}s^\prime \right).
	\end{align}
	We assume $A(x)>0$ for $x>0$. \eqref{eq:2.15} is bounded from above by
	\begin{align}
		\sum_{j=0}^{t-1} \int_j^{j+1}\P\left( \exists k\leq n(t) : x_k(\lceil s^\prime \rceil )<U_{A,\g/2}(\lceil s^\prime\rceil) \bigg| \mathcal{T}=s^\prime\right) \P\left( \mathcal{T}\in \mathrm{d}s^\prime \right).
	\end{align}
	It remains to show that 
	\begin{align}
		\P\left( \exists k\leq n(t) : x_k(\lceil s^\prime \rceil )<U_{A,\g/2}(\lceil s^\prime \rceil) \bigg| \mathcal{T}=s^\prime\right)
	\end{align}
	tends to $0$ uniformly in $s^\prime$, as $C\rightarrow \infty$. By the definition of $\mathcal{T}$, this probability is bounded from above by the probability that the predecessor at time $\lceil s^\prime \rceil$ of the maximum at time $s^\prime$ makes a downward jump smaller than
	\begin{align}
		&U_{A,\g}(\lceil s^\prime \rceil )-U_{A,\g /2}( s^\prime )<\sqrt{2}t\left(A(\lceil s^\prime \rceil/t) -A(s^\prime/t)\right)\\ &+\big((t A(\lceil s^\prime \rceil /t) \wedge (t-A(\lceil s^\prime \rceil /t)t))\vee C\big)^{\g}-\big((t A(\lceil s^\prime \rceil /t) \wedge (t-A(\lceil s^\prime \rceil /t)t))\vee C\big)^{\g/2} \nonumber,
	\end{align}
	whose probability is, by the Markov property of BBM and Markov's inequality, bounded from above by
	\begin{align}\label{eq:2.19}
		\eee^{\lceil s^\prime \rceil - s^\prime }\P\left(x_1(\lceil s^\prime \rceil ) -x_1( s^\prime )<  U_{A,\g}(\lceil s^\prime \rceil )- U_{A,\g/2}( s^\prime ) \right)  .
	\end{align}
	Now, $\eee^{\lceil s^\prime \rceil - s^\prime }<e$ and we also have $A( s^\prime /t)>0$.
	Note further that $x_1 \left( \lceil s^\prime \rceil\right)-x_1\left( s^\prime \right)\sim \mathcal{N}\left(0,t \left[ A\left( \lceil s^\prime\rceil/t \right) -A\left(s^\prime/t\right)\right] \right)$ and that, for $s^\prime$ so that $0< t A(s^\prime/t)\leq t/2$,
	\begin{align}
		&U_{A,\g}(\lceil s^\prime \rceil )-m_A(s^\prime)-U_{A,\g /2}( s^\prime ) +m_A(\lceil s^\prime \rceil)\nonumber\\
		&\qquad =-\left((t A(\lceil s^\prime \rceil /t)) \vee C\right)^{\g}\left[ 1- \frac{\left((t A(\lceil s^\prime \rceil/t ))\vee C\right)^{\g/2}}{\left((t A( s^\prime/t))\vee C \right)^{2}}  \right]\leq -\frac{1}{2}C^\g,
	\end{align}
	where we choose $C>0$ sufficiently large to get the factor $1/2$. Similarly, for $s^\prime$ with $t/2< tA(s^\prime/t)\leq t,$
	\begin{align}
		&U_{A,\g}(\lceil s^\prime \rceil )-m_A(s^\prime)-U_{A,\g /2}( s^\prime ) +m_A(\lceil s^\prime \rceil)
		=-\left((t-tA(s^\prime/t))\vee C\right)^{\g}\nonumber\\& \qquad \times\left( 1- \frac{\left((t- tA(\lceil s^\prime \rceil/t))\vee C \right)^{\g/2}}{\left((t-tA(s^\prime /t))\vee C\right)^{\g}}\right) \leq -\frac{1}{2}C^{\g}.
	\end{align}
	Therefore \eqref{eq:2.19} is bounded from above by
	\begin{align}
		\eee^{\lceil s^\prime \rceil -s^\prime} 
		\mathbb{P}\left( X \leq \sqrt{2t} \left(\sqrt{ \lceil s^\prime \rceil A(\lceil s^\prime \rceil /t)}-\sqrt{s^\prime A(s^\prime /t)} -\frac{1}{2}C^{\g} \right)\right),
	\end{align}
	where $X\sim \mathcal{N}\left(0, t\left(A(\lceil s^\prime \rceil /t)-A(s^\prime /t) \right)\right)$. Thus, this tends to $0$ as $C\rightarrow \infty$, uniformly in $s^\prime$, which concludes the proof.
\end{proof}
\section{Phase B1: Fluctuation of the partition function}

In this section, we analyze the behavior of the partition function in Phase $B_1$. The main result in phase $B_1$ is Theorem \ref{paseB1}. 
\subsection{Proof of Theorem \ref{paseB1}}
A first step  in the proof of Theorem \ref{paseB1} is the following proposition.

\begin{proposition}\label{Prop.B1.1}
	For any $\delta>0$ and $\e>0$, there exists $r_0$ such that for all $r\geq r_0$ and all $t$ sufficiently large
	\begin{equation}
		\P\left(\left\vert  \XX_{\beta,\rho}(t) - \E\left(\XX_{\beta,\rho}(t) | \FF_{r}\right)\right\vert \geq \delta \eee^{t\left(1+\frac{\s^2-\t^2}{2}\right)}\right)\leq \epsilon.
	\end{equation}
\end{proposition}
We postpone the proof of Proposition \ref{Prop.B1.1} to the end of the section.
\begin{remark}
	\label{rem:decomposition-dependent-and-independent-components}
	It holds that
	\begin{align}\Eq(cor.1)
		Y(t)\overset{\mathrm{D}}{=} \rho X(t)  +\sqrt{1-\rho^2} Z(t)
		,
	\end{align}
	where ``$\overset{\mathrm{D}}{=}$'' denotes equality in distribution and $Z(t)
	:= (z_i(t))_{i\leq n(t)}$ is a CREM independent of $X(t)$
	and defined with respect to the same underlying GW process. This is used repeatedly throughout this article to handle the correlation between $X$ and $Y$.
\end{remark}
Next, we show that the conditional expectation of $\XX_{\beta,\rho}(t)$ on $\FF_r$ is close to
\begin{equation}
	\eee^{t\left(1+2i\rho\s\t+\frac{\s^2-\t^2}{2}\right)}\MM_{\s_b\s,\s_b\t}(r).
\end{equation}

\begin{lemma}\label{Lem.B1.2}
	Let  $(\sigma_b \tilde x_k(r))_{k\leq n(r)}$ and $(x_k(r))_{k\leq n(r)}$ be coupled as described in Remark~\ref{rem:coupling} and let $\widehat\MM_{\s_b\s,\s_b\t}(r)$ be the corresponding coupled realization of $\MM_{\s_b\s,\s_b\t}(r)$ . For any $\delta>0$ and $\e>0$, there exists $r_0$ such that for all $r\geq r_0$ and all $t$ sufficiently large
	\be\Eq(B1.6)
	\P\left( \left|\eee^{-t\left(1 +\frac{\s^2-\t^2}{2}\right)}\E\left(\XX_{\beta,\rho}(t) \mid \FF_{r}\right)-\widehat\MM_{\s_b\s,\s_b\t}(r)\right|>\delta \right)<\e.
	\ee
\end{lemma}

Now, modulo the proof of Proposition \ref{Prop.B1.1} and Lemma \ref{Lem.B1.2} we are in the position to prove Theorem \ref{paseB1}.
\begin{proof}[Proof of Theorem \ref{paseB1}]
	This follows immediately from  Proposition \ref{Prop.B1.1} and Lemma \ref{Lem.B1.2}.  
\end{proof}
To only prove convergence in distribution, one might use the following lemma
instead of Lemma \ref{Lem.B1.2}.  
\begin{lemma}\label{Lem.B1.2dist} Let $\beta=\s+i\t\in B_1$. Then
	$\eee^{-r-tA(r/t)(\s^2-\t^2)/2}\E\left(\XX_{\beta,\rho}(t) \mid \FF_r\right)$ converges in distribution to $\MM_{\s_b\s,\s_b\t}(r)$, as $t\to\infty$.
\end{lemma}
\begin{proof}
	Note that the covariance function of $(x_k(r))_{k\leq n(r)}$ converges as $t\uparrow \infty$ to the one of $(\sigma_b \tilde x_k(r))_{k\leq n(r)}$, where $(\tilde x_k(r))_{k\leq n(r)}$ are distributed as particles of a standard branching Brownian motion. The claim then follows from the continuous mapping theorem. 
\end{proof}
A direct consequence of Proposition \ref{Prop.B1.1} and Lemma \ref{Lem.B1.2dist} is the following.
\begin{corollary}\label{paseB1.1}
	Let $\beta =\s+i\t$ with $\beta\in B_1, $, and $\rho\in[-1,1]$. Then,
	\be
	\eee^{-t\left(1+2i\rho\s\t+\frac{\s^2-\t^2}{2}\right)}\XX_{\b,\rho}(t) \Rightarrow \mathcal{M}_{\s_b\s,\s_b\t},\qquad \text{as }t\uparrow \infty,
	\ee
	where convergence is in law.
\end{corollary}

\subsection{Proof of Proposition \ref{Prop.B1.1} and Lemma \ref{Lem.B1.2}}

We  first prove Lemma \ref{Lem.B1.2}.
\begin{proof}[Proof of Lemma \ref{Lem.B1.2}]
	The probability in \eqref{B1.6} is bounded from above by
	\begin{eqnarray}\label{almost.mart}
		&&\E\Big(\eee^{-2r}\sum_{k,k'\leq n(r)}
		\big(\eee^{-(\s^2-\t^2)tA(r/t)}e^{\s x_k(r) +i\tau y_k(r)+\s x_{k'}(r)+i\tau y_{k'}(r)} \nonumber\\
		&&\qquad\qquad\qquad- \eee^{-\frac{\s^2-\t^2}{2}(tA(r/t)+\s_b^2 r)}\eee^{\s \s_b\tilde x_k(r)+i\t \s_b\tilde y_{k}(r)+\s x_{k'}(r)-i\t y_{k'}(r)}
		\nonumber\\
		&&\qquad\qquad\qquad
		-\eee^{ -\frac{\s^2-\t^2}{2}(tA(r/t)+\s_b^2 r)}\eee^{\s \s_b\tilde x_{k'}(r)-i\t \s_b\tilde y_{k'}(r)+\s x_{k}(r)+i\t y_{k}(r)}\nonumber\\
		&&\qquad\qquad\qquad+ \eee^{ -(\s^2-\t^2)\s_b^2r}
		\eee^{\s\s_b\tilde x_k(r)+\t\s_b\tilde y_k(r)+\s\s_b \tilde x_{k'}(r)-i\t\s\s_b \tilde y_{k'}(r)}\big)\Big).
	\end{eqnarray}
	First, note that the two summands with negative signs are complex conjugates of each other and thus their sum is real valued and the other terms are also real numbers. Hence, we do not have to worry about the imaginary parts in what follows, as they cancel out.
	We show that the expectation of the four summands is equal. Second, note that since $A$ is twice differentiable in a neighborhood of zero, by a first order Taylor expansion
	\begin{equation}\label{mar.6}
		A(s/t)=\s_b^2 s+ O(s^2/t).
	\end{equation}
	This implies that the multiplicative prefactors of the four summands are all
	\begin{equation}
		\eee^{-(\s^2-\t^2)\s_b^2r+O(r^2/t)} \to \eee^{ -(\s^2-\t^2)\s_b^2r}, \quad \text{as } t \to \infty.
	\end{equation}
	Next, we compute
	\begin{equation}\label{mar.1}
		\E\Big(\sum_{k,k'\leq n(r)}
		\eee^{\s \s_b\tilde x_k(r)+i\t \s_b\tilde y_{k}(r)+\s x_{k'}(r)-i\t y_{k'}(r)}
		\Big).
	\end{equation}
	As the correlation coefficient between the real part and the imaginary part is $\rho$, using representation \eqref{cor.1}, the expectation in \eqref{mar.1} is equal to 
	\begin{equation}
		\E\Big(\sum_{k,k'\leq n(r)}
		\eee^{(\s+i\rho\t) \s_b\tilde x_k(r)+i\sqrt{1-\rho^2}\t  \s_b\tilde z_{k}(r)+(\s-i\rho \t) x_{k'}(r)-i\sqrt{1-\rho^2}\t z_{k'}(r)}
		\Big).
	\end{equation}
	Let $d\equiv d(x_k.x_{k'})$.  Then, conditional on the underlying Galton-Watson tree,
	\begin{multline}\label{mar.2}
		\E\left(\eee^{i\sqrt{1-\rho^2}\t  \s_b\tilde z_{k}(r))-i\sqrt{1-\rho^2}\t z_{k'}(r)}\right)=\E\left(e^{i\sqrt{1-\rho^2}\t(\sqrt{tA(d/t)}-\s_b\sqrt{d})z_1}\right)
		\\
		\times
		\E\left(\eee^{i\sqrt{1-\rho^2}\t \sqrt{tA(r/t)-tA(d/t)} z_2}\right)\E\left(\eee^{-i\sqrt{1-\rho^2}\t \s_b \sqrt{r-d} z_3}\right),
	\end{multline}
	where $z_1,z_2,z_3$ are standard Gaussians.
	\eqref{mar.2} is equal to
	\begin{equation}\label{mar.4}
		\eee^{-(1-\rho^2)\frac{\t^2}{2}(\sqrt{tA(d/t)}-\s_b\sqrt{d})^2}
		\eee^{-(1-\rho^2)\frac{\t^2}{2} (tA(r/t)-tA(d/t)+ \s_b^2(r-d))}.
	\end{equation}
	By \eqref{mar.6}, this coincides (in absolute value)  with
	\begin{equation}\label{mar.7}
		\E\left(\eee^{i\sqrt{1-\rho^2}\t  \s_b\tilde z_{k}(r))-i\sqrt{1-\rho^2}\t \s_b\tilde z_{k'}(r)}\right).
	\end{equation}
	Next, we compute the expectation over $X$ and $\tilde X$, again conditional on the underlying tree.
	\begin{multline}
		\label{mar.3}
		\E\Big(
		\eee^{(\s+i\rho\t) \s_b\tilde x_k(r)+(\s-i\rho \t) x_{k'}(r)}\Big)= 
		\E\left(\eee^{((\s+i\rho\t)\sqrt{tA(d/t)}+
			(\s-i\rho\t)\s_b\sqrt{d})x_1}\right)
		\\
		\times
		\E\left(\eee^{(\s+i\rho\t) \sqrt{tA(r/t)-tA(d/t)} x_2}\right)\E\left(\eee^{(\s-i\rho\t) \s_b \sqrt{r-d} x_3}\right),
	\end{multline}
	where $x_1,x_2,x_3$ are standard Gaussians. 
	This is equal to
	\begin{equation}\label{mar.5}
		\eee^{\frac{1}{2}\left((\s+i\rho \t)\sqrt{tA(d/t)}+ (\s-i\rho\t)\s_b\sqrt{d}\right)^2}
		\eee^{\frac{(\s+i\rho \t)^2}{2}(tA(r/t)-tA(d/r))}
		\eee^{\frac{(\s-i\rho \t)^2}{2}\s_b^2(r-d)}.
	\end{equation}
	Again by \eqref{mar.6}, this coincides in absolute value with
	\begin{equation}\label{mar.8}
		\eee^{O(r^2/t)}\E\Big(
		\eee^{(\s+i\rho\t) \s_b\tilde x_k(r)+(\s-i\rho \t)\s_b\tilde x_{k'}(r)}
		\Big).
	\end{equation}
	\eqref{mar.7} and \eqref{mar.8} now imply the desired result together with the observation made above that we do not need to take care of the imaginary parts, as $\eee^{O(r^2/t)}$ converges to one as $t$ tends to infinity.
	
	
	
\end{proof}

\begin{proof}[Proof of Proposition \ref{Prop.B1.1}]
		Note that it suffices to bound 
		\be\Eq(B1.9)
		\P\left(\left|  \sum_{k=1}^{n(t)} \eee^{\s x_k(t)+i\t y_k(t)}\1_{ x_k\in \UU_{t,\g} } - \E\left(\sum_{k=1}^{n(t)} \eee^{\s x_k(t)+i\t y_k(t)}\1_{ x_k\in \UU_{t,\g} } \mid \FF_{r}\right)\right|\geq \delta \eee^{t\left(1+\frac{\s^2-\t^2}{2}\right)}\right),
		\ee
		by Lemma~\ref{Lem.B1.1} and 
		\begin{lemma}\label{Lem.expec} Let $\e>0$ and $r,t\in \R_+$ satisfying $r\ll t$. Then, there exists $C(r,\e)>0$ such that 
			\begin{align}
				&\eee^{-t\left(1+\frac{\s^2-2i\t \s+
						\t^2}{2}\right)}\E\left(\sum_{k=1}^{n(t)} \eee^{\s x_k(t)+i\t y_k(t)}\1_{ x_k\in \UU_{t,\g} } \mid \FF_{r}\right)\nonumber\\
				&\geq
				(1-\e) \eee^{-t\left(1+\frac{\s^2-2i\t \s+
						\t^2}{2}\right)}
				\E\left(\sum_{k=1}^{n(t)} \eee^{\s x_k(t)+i\t y_k(t)}  \mid \FF_{r}\right)
			\end{align}
			with probability $>1-\e$.
		\end{lemma}
		\begin{proof}[Proof of Lemma~\ref{Lem.expec}]
			First, note that
			\bea \label{expec.1}
			&&\eee^{-t\left(1+\frac{\s^2+2i\s\t -\t^2}{2}\right)}\E\left(\sum_{k=1}^{n(t)} \eee^{\s x_k(t)+i\t y_k(t)} \mid \FF_{r}\right)
			\nonumber\\
			&&
			= 
			\eee^{-tA(r/t) \left(1+\frac{\s^2+2i\s\t -\t^2}{2}\right)} \sum_{k=1}^{n(r)} \eee^{\s x_k(tA(r/t))+i\t y_k(tA(r/t))}.
			\eea
			Next, we compute  $\E\left(\sum_{k=1}^{n(t)} \eee^{\s x_k(t)+i\t y_k(t)}\1_{ x_k\in \UU_{t,\g} } \mid \FF_{r}\right)$. By the branching property and the many-to-one lemma, the expectation is equal to
			\begin{eqnarray}\label{expec.3}
				&&\sum_{k=1}^{n(r)} e^{t-r} \eee^{\s x_k(tA(r/t))+i\t y_k(tA(r/t)) }
				\1_{\left\{x_k(r')\leq \sqrt{2r'tA(r'/t)} +\big(C\vee( tA(r'/t) \land (t-tA(r'/t)))\big)^{\g},\forall r'\leq r\right\}}\nonumber\\
				&\times&\E\left( \eee^{\s x(t-r)+i\t y(t-r)}\1_{\left\{x(s-r)+x_k(r) \leq \sqrt{2 s tA(s/t)}+ \big(C\vee( tA(s/t) \land (t-tA(s/t)))\big)^{\g},\, \forall s>r\right\}}\right)
			\end{eqnarray}
			where $(x(s-r))_{r<s<t}$ and  $(y(s-r))_{r<s<t}$ are time-inhomogeneous Brownian motions with variance $\hat s\equiv tA(s/t)-tA(r/t)$ starting from $0$ with correlation coefficient $\rho$. It holds that $(x(s-r))_{r<s<t}\stackrel{d}{=} (B(\hat s))_{0\leq \hat s\leq t-tA(r/t)}$, where $B$ is a standard Brownian motion.
			Hence,
			\begin{equation}
				\xi(\hat s)= B(\hat s)-\frac{\hat s}{t-tA(r/t)} B(t-tA(r/t)),
			\end{equation}
			is a Brownian bridge from $0$ to $0$ in time $t-tA(r/t)$ which is independent of $B(t-tA(r/t))$.
			We rewrite
			\begin{align}\label{expec.2}
				x(s-r)+x_k(r)=\xi(\hat s) + \frac{t-tA(s/t)}{t-tA(r/t)}x_k(r) + \frac{\hat s}{t-tA(r/t)} (B(t-tA(r/t))+x_k(r))
			.
		\end{align}
		Set $ f_\g(C,s):=\big(C \vee (tA(s/t)\land (t-tA(s/t))\big)^\g$
		The requirement that $x_k\in \mathcal{U}_{t,\g}$ implies that $\{x(r)+x(t-r)\leq \sqrt{2}t\}$ for $t$ sufficiently large. In our new notation, this can be rewritten as $\{B(t-tA(r/t))+x_k(r)\leq \sqrt{2}t\}$. Thus, the event 
		\begin{eqnarray}\label{expec.5}
			\left\{\xi(\hat s) + \frac{t-tA(s/t)}{t-tA(r/t)}x_k(r) + \frac{\sqrt{2} t\hat s}{t-tA(r/t)} 
			\leq \sqrt{2 s tA(s/t)}+  f_\g(C,s) \,\forall s>r
			\right\} 
		\end{eqnarray}
		is a subset of the event in the second indicator function in \eqref{expec.3}.  As $x_k\in \mathcal{U}_{t,\g}$, implies $ x_k(r)\leq \sqrt{2trA(r/t)} + (C \lor rA(r/t))^{\g}$, the event in \eqref{expec.5} contains 
		\begin{eqnarray} \label{expec.5.1}
			\Big\{\xi(\hat s) 
			\leq \sqrt{2  tsA(s/t)}- \frac{\sqrt{2}\hat s}{1-A(r/t)}+  f_\g(C,s)\nonumber\\
			-\sqrt{2rtA(r/t)}-(C \lor rA(r/t))^{\g}\,,\forall s>r
			\Big\},
		\end{eqnarray}
		as a subset.
		As $tA(s/t)\leq s$ and $\sqrt{2 s tA(s/t)}-\sqrt{2} \frac{\hat s}{1-A(r/t)}\leq O(tA(r/t))$, it follows that 
		\begin{equation}\label{expec.5.2}
			\mathcal{L}_r=\left\{\xi(\hat s) 
			\leq  f_\g(C,s) + O\left(\sqrt{rtA(r/t)}\right)\forall s>r
			\right\}
		\end{equation}
		is a subset of \eqref{expec.5.1}
		The expectation in \eqref{expec.2} (on the event $\mathcal{L}_r$) is bounded from below by
		\begin{eqnarray}\label{expec.6}
			\E\left( \eee^{\s x(t-r)+i\t y(t-r)}\right)
			\P\left(\xi(\hat s) 
			\leq  f_\g(C,s) + O\left(\sqrt{rtA(r/t)}\right)\forall s>r\right).
		\end{eqnarray}
		If we choose $C=C(r)$ large enough, then we can bound the probability in \eqref{expec.6} from above by
		\begin{equation}
			\P\left( \xi(\hat s) 
			\leq \big(C \vee (\hat s \land (t-A(r/t)-\hat s))\big)^\g, \forall 0\leq \hat s \leq t-tA(r/t)\right)\geq 1-\epsilon.
		\end{equation}
		Note that by Lemma \ref{Lem.B1.1} the probability of the event
		\be
		\left\{x_k(r')\leq \sqrt{2r'tA(r'/t)} +\big(C\vee( tA(r'/t) \land (t-tA(r'/t)))\big)^{\g},\forall r'\leq r,k\leq n(r)\right\}
		\ee
		is at least $1-\e$, for $C$ large enough. On this event, \eqref{expec.3} is bounded from below by $(1-\e)$ times the expression in \eqref{expec.1}. This implies the claim of Lemma~\ref{Lem.expec}.
	\end{proof}
	
	We continue the proof of Proposition \ref{Prop.B1.1}. By Chebyshev's inequality, the probability in \eqv(B1.9) is bounded from above by $\delta^{-2} \eee^{-t\left(2+\s^2-\t^2\right)} $ times
	\be\Eq(B1.10)
	\E\left(\left|  \sum_{k=1}^{n(t)} \eee^{\s x_k(t)+i\t y_k(t)}\1_{ x_k\in \UU_{t,\g} } - \E\left(\sum_{k=1}^{n(t)} \eee^{\s x_k(t)+i\t y_k(t)}\1_{ x_k\in \UU_{t,\g} } \mid \FF_{r}\right)\right|^2\right).
	\ee
	To continue estimating \eqref{B1.10}, we need the indicators
	\begin{equation}
		\label{B1.101}
		I_{k,j} := \1\{ x_k\left(r\right) + x^k_j\left(s-r\right) \leq U_{A,\g}(s), \forall s \in [r,t] \},
	\end{equation}
	for $k \leq n(r)$ and $j \leq n^k(t-r)$,
	where $n^k(t-r)$ denotes the number of offspring of the particle $x_k(r)$ alive at time $t$ and $x^k_j(s)$ are the positions of the $j^\mathrm{th}$ offspring of the particle $x_k(r)$ alive at time $r+s$. Due to the branching property, $n^k(s)$, for $k\leq n(r)$, the total number of offspring alive at time $r+s$ of particles $x_k(r)$, are independent for different $k$.
	Noting that $|x|^2=x\overline{x}$ for any $x\in \C$ and using representation \eqref{cor.1}, we rewrite the expectation in \eqv(B1.10) as
	\bea\Eq(B1.11)
	&&\E\Biggl(\sum_{k,k'=1}^{n(r)} \eee^{(\s+i\t\rho)x_k(r)+i\t\sqrt{1-\rho^2}z_k(r)+(\s-i\t \rho)x_{k'}(r)-i\t \sqrt{1-\rho^2}z_{k'}(r)}\nonumber\\ && \quad\times\Biggl(
	\E\left(\sum_{j=1}^{n^{k}(t-r)}\sum_{j'=1}^{n^{k'}(t-r)} \eee^{(\s+i\t\rho)x^k_j(t-r)+i\t\sqrt{1-\rho^2}z^k_j(t-r)}I_{k, j}\right.\nonumber\\ &&
	\left.\qquad \qquad \quad\times \eee^{(\s-i\t \rho)x^{k'}_{j'}(t-r)-i\t \sqrt{1-\rho^2}z^{k'}_{j'}(t-r)-2t-(\s^2-\t^2)t}
	I_{k^\prime, j^\prime} \right)\nonumber\\
	&&\quad-\left|\E\left(\sum_{j=1}^{n^{k}(t-r)}\eee^{(\s+i\t\rho)x_j^k(t-r)+i\t\sqrt{1-\rho^2}z_k(t-r)-t-(\s^2-\t^2)t/2} I_{k,j} \right)\right|^2\Biggr)\Biggr).
	\eea
	We set $\lambda=(\s+i\t \rho)$ and observe that for particles with $d(x_k,x_k')\leq r$, $x^k$ and $x^{k'}$ as well as $z^k$ and $z^{k'}$ are independent and thus, for such particles,
	\begin{multline}\label{nico.1}
		\E\left(\sum_{j=1}^{n^{k}(t-r)}\sum_{j'=1}^{n^{k'}(t-r)} \eee^{\lambda x^k_j(t-r)+i\t\sqrt{1-\rho^2}z^k_j(t-r)+(\s-i\t \rho)}I_{k, j}\right.
		\\ \left. \times 
		\eee^{\overline{\lambda}x^{k'}_{j'}(t-r)-i\t \sqrt{1-\rho^2}z^{k'}_{j'}(t-r)-2t-(\s^2-\t^2)t}
		I_{k^\prime, j^\prime} \right)\\
		=
		\left|\E\left(\sum_{j=1}^{n^{k}(t-r)}\eee^{\lambda x_j^k(t-r)+i\t\sqrt{1-\rho^2}z_k(t-r)-t-(\s^2-\t^2)t/2} I_{k,j} \right)\right|^2.
	\end{multline}
	By this observation, \eqref{B1.11} is equal to
	\bea\Eq(B1.112)
	&&\E\Biggl(\sum_{k=1}^{n(r)} \eee^{2\s x_k(r)} \Biggl(
	\E\left(\sum_{j=1}^{n^{k}(t-r)}\sum_{j'=1}^{n^{k'}(t-r)} \eee^{\lambda x^k_j(t-r)+i\t\sqrt{1-\rho^2}z^k_j(t-r)} I_{k, j}\right.\nonumber\\ && \left. \qquad\times \eee^{(\s-i\t \rho)x^{k'}_{j'}(t-r)-i\t \sqrt{1-\rho^2}z^{k'}_{j'}(t-r)-2t-(\s^2-\t^2)t}
	I_{k^\prime, j^\prime} \right)\nonumber\\
	&&\quad-\left|\E\left(\sum_{j=1}^{n^{k}(t-r)}\eee^{\lambda x_j^k(t-r)+i\t\sqrt{1-\rho^2}z_k(t-r)-t-(\s^2-\t^2)t/2} I_{k,j} \right)\right|^2\Biggr)\Biggr).
	\eea
	This is obviously bounded from above by
	\bea\label{nico.4}
	&&\E\left(\sum_{k=1}^{n(r)} \eee^{\lambda x_k(r)+i\t\sqrt{1-\rho^2}z_k(r)+\lambda x_{k'}(r)-i\t \sqrt{1-\rho^2}z_{k'}(r)} \E\left(\sum_{j,j'=1}^{n^{k}(t-r)}e^{\lambda x_j(t-r)+i\t\sqrt{1-\rho^2}z_j(t-r)}\right.\right. \nonumber \\ && \left. \left.
	\quad \times
	\eee^{\overline{\lambda}x_{j'}^k(t-r)-i\t \sqrt{1-\rho^2}z_{j'}^k(t-r)-2t-(\s^2-\t^2)t}I_{k,j} I_{k,j^\prime}\right) \right) .
	\eea
	We observe that  \eqref{nico.4} can be rewritten as a second moment with an additional indicator function that the particles branched after time r. Namely, \eqref{nico.4} is equal to $\eee^{-2t-(\s^2-\t^2)t}$ times
	\be\Eq(B1.12)
	\E\left(\sum_{k,k'=1}^{n(t)} \eee^{\l x_k(t)+\bar\l x_{k'}(t)+i\t \sqrt{1-\rho^2}(z_k(t)-z_{k'}(t))}\1_{x_k,x_{k'}\in  \UU_{t,\g} }\1_{d(x_k,x_{k'})\geq r}\right). 
	\ee
	By keeping the constraint on the height of the particles only at $d(x_k,x_{k'})$, we bound \eqv(B1.12) from above by
	\be\Eq(B1.13)
	\E\left(\sum_{k,k'=1}^{n(t)} \eee^{\lambda x_k(t)+\overline{\lambda} x_{k'}(t)+i\t \sqrt{1-\rho^2}(z_k(t)-z_{k'}(t))}\1_{x_k\left(d(x_k,x_{k'})\right)\leq U_{A,\g}\left(d(x_k,x_{k'})\right) }\1_{d(x_k,x_{k'})\geq r}\right). 
	\ee
	As $Z$ is independent from $X$, \eqv(B1.13) is equal to
	\begin{align}\label{eq:extra.2}
		\E\left(\sum_{k,k'=1}^{n(t)}  \eee^{-\t^2(1-\rho^2)(t-tA(d(x_k,x_{k'})/t))}\eee^{\lambda x_k(t)+\overline{\lambda} x_{k'}(t)}\1_{x_k\left(d(x_k,x_{k'})\right)\leq U_{A,\g}\left(d(x_k,x_{k'})\right) }\1_{d(x_k,x_{k'})\geq r}\right). 
	\end{align}
	By the many-to-two lemma and \eqref{eq:extra.2},  \eqv(nico.4) is  bounded from above by
	\bea\Eq(B1.14)
	&&\int_{r}^t \dd s \eee^{2t-s-2t-(\s^2-\t^2)t -\t^2(1-\rho^2)(t-tA(s(/t))} \int_{-\infty}^{U_{A,\g}(s)}\frac{\dd x}{\sqrt{2\pi tA(s/t)}} \eee^{2\s x-\frac{x^2}{2 tA(s/t)}}\nonumber\\
	&&\times \int_{-\infty}^{\infty}\frac{\dd y \eee^{\l y-\frac{y^2}{2t\left(1-A(s/t)\right)}}}{\sqrt{2\pi t\left(1-A(s/t)\right)}} 
	\int_{-\infty}^{\infty}\frac{\dd y' \eee^{\overline{\l} y'-\frac{(y')^2}{2t\left(1-A(s/t)\right)}}}{\sqrt{2\pi t\left(1-A(s/t)\right)}}.
	\eea
	Using that the integrals with respect to $y,y'$ are Fourier transforms of Gaussian random variables, we get 
	\bea\Eq(B1.15)
	\int_{-\infty}^{\infty}\frac{\dd y \eee^{\l y-\frac{y^2}{2t\left(1-A(s/t)\right)}}}{\sqrt{2\pi t\left(1-A(s/t)\right)}} 
	\int_{-\infty}^{\infty}\frac{\dd y' \eee^{\overline{\l} y'-\frac{(y')^2}{2t\left(1-A(s/t)\right)}}}{\sqrt{2\pi t\left(1-A(s/t)\right)}} 
	=\eee^{\left(\s^2-\rho^2\t^2\right)t\left(1-A(s/t)\right)}.
	\eea
	Hence, \eqv(B1.14) is equal to
	\bea\Eq(B1.16)
	&&\int_{r}^t \dd s \eee^{-s-(\s^2-\t^2)tA(s/t)} \int_{-\infty}^{U_{A,\g}(s)}\frac{\dd x}{\sqrt{2\pi tA(s/t)}} \eee^{2\s x-\frac{x^2}{2 tA(s/t)}}\nonumber\\
	&=&\int_{r}^t \dd s \eee^{-s-(\s^2-\t^2) tA(s/t)+2\s^2 tA(s/t)} \int_{-\infty}^{U_{A,\g}(s)}\frac{\dd x}{\sqrt{2\pi tA(s/t)}} \eee^{-\frac{\left(x-2\s t A(s/t)\right)^2}{2 tA(s/t)}}\nonumber\\
	&=& \int_{r}^t \dd s  \eee^{-s+(\s^2+\t^2)tA(s/t)} \int_{-\infty}^{U_{A,\g}(s)-2\s t A(s/t)}\frac{\dd z}{\sqrt{2\pi tA(s/t)}} \eee^{-\frac{z^2}{2 tA(s/t)}}.
	\eea
	Next, we distinguish two cases. If $U_{A,\g}(s)-2\s t A(s/t)\geq 0$, we bound the integral with respect to $z$ by one. Otherwise, we use a Gaussian tail bound, to bound \eqv(B1.16) from above by
	\bea\Eq(B1.17)
	&&\int_{r}^t \dd s \eee^{-s+(\s^2+\t^2)tA(s/t)}\1_{U_{A,\g}(s)-2\s t A(s/t)\geq 0}\nonumber\\
	&&+ \int_{r}^t \dd s \eee^{-s+(\s^2+\t^2) tA(s/t)}\1_{U_{A,\g}(s)-2\s t A(s/t)< 0} \eee^{-\frac{\left(U_{A,\g}(s)-2\s t A(s/t)\right)^2}{2tA(s/t)}}.
	\eea
	We start by bounding the first summand in \eqv(B1.17). We distinguish two cases. For $\t^2+\s^2\leq 1$ the first summand is bounded from above by 
	\be\Eq(B1.19)
	\int_{r}^t \dd s \eee^{-\epsilon s}\leq \frac{1}{\epsilon}\eee^{-\epsilon r},
	\ee
	for some $\e>0$, where we used Assumption~\ref{weak}. Otherwise, observe that the event in the indicator of the first summand implies  $s\geq 2\s^2 tA(s/t)+ O\left((tA(s/t))^\g\right)$.
	Hence, on this event
	\be\Eq(B1.18)
	-s+(\s^2+\t^2) tA(s/t)\leq (\t^2-\s^2) tA(s/t).
	\ee
	Now, as for $\b\in B_1$ and $\s^2+\t^2>1$, we have $\s>\t$,  in this case the first summand in \eqv(B1.17) is bounded from above by
	\be\Eq(B1.19.2)
	\int_{r}^t \dd s \eee^{-\epsilon t A(s/t)}\leq C\eee^{-\epsilon tA(r/t)},
	\ee
	for some $\epsilon>0$.
	Using \eqv(B1.18), there exists $\epsilon>0$ such that the first summand in \eqv(B1.17) is bounded from above by \eqref{B1.19}.
	The second summand in \eqv(B1.17) is bounded from above by
	\begin{eqnarray}\Eq(B1.20)
		\int_{r}^t \dd s \eee^{-2s-(\s^2-\t^2)tA(s/t)+2\s\sqrt{2stA(s/t)}+O\left(\sqrt{s}\lor (tA(s/t))^\gamma\right)}\nonumber\\
		= \int_{r}^t \dd s e^{\t^2 tA(s/t)}\eee^{-\left(\sqrt{2s}-\s\sqrt{tA(s/t)}\right)^2+O\left(\sqrt{s}\lor (tA(s/t))^\gamma\right)}.
	\end{eqnarray}
	As $\s+\t <\sqrt{2}$ and $\sqrt{tA(s/t)}<\sqrt{ s }$, the exponent in \eqv(B1.20) is bounded from above by $-\epsilon' s + O(\sqrt{ s })$, for some $\epsilon'>0$ depending on $\s$. Hence, \eqv(B1.20) is bounded from above by
	\be\Eq(B1.21)
	\int_{r}^t \dd s \eee^{-\e' s}\leq \frac{1}{\e'}\eee^{-\e'r}.
	\ee
	Combining \eqv(B1.19), \eqv(B1.19) and \eqv(B1.21), we get the claim of the lemma.
\end{proof}

\section{Phase B2}
In phase $B_2$, as already mentioned in the introduction, the behaviour of the partition function is determined by the extrema of $X$. This is made precise in Proposition \ref{Prop.B2.1}. Theorem \ref{phaseB2.2} and \ref{phaseB2.1} are then proven using the continuous mapping theorem.
\begin{proposition}\label{Prop.B2.1}
	If $\rho \in [-1,1]$ and $\beta \in B_2$, then for all $\delta, \epsilon>0$ there is $B_0$ such that for all $B\geq B_0$ and all $t$ large enough,
	\begin{align}\label{B2.016}
		\P\left( \left|\sum_{k=1}^{n(t)} \eee^{\s x_k(t)+i\t y_k(t)-\s m(t)}\1_{x_k(t)<m(t)-B}\right|>\d\right)<\epsilon.
	\end{align}
\end{proposition}
\begin{proof}[Proof of Proposition~\ref{Prop.B2.1}:]
	We distinguish two cases:
	\begin{equation}
		\mbox{Case (a):}\, \b\in B_2, \s>\sqrt{2},\quad \mbox{Case (b):}\,\b\in B_2, \s\leq\sqrt{2}.
	\end{equation}
	We first prove Proposition \ref{Prop.B2.1} in Case (a). By Markov's inequality and the triangle inequality, the probability in \eqref{B2.016} is bounded from above by
	\begin{equation}\label{B2.017}
		\d^{-1}\E\left( \sum_{k=1}^{n(t)} \eee^{\s (x_k(t)-m(t))}\1_{x_k(t)<m(t)-B}   \right).
	\end{equation}
	Now, the expectation in \eqref{B2.017} is equal to
	\begin{equation}\label{B2.018}
		\eee^{t-\s m(t)}\int_{-\infty}^{m(t)-B}\eee^{\s x-\frac{x^2}{2t}}\frac{\dd x}{\sqrt{2\pi t}}
		=  \eee^{t-\s m(t)+\sigma^2 t/2}\int_{-\infty}^{m(t)-B}\eee^{-\frac{(x-\s t)^2}{2t}}\frac{\dd x}{\sqrt{2\pi t}}.
	\end{equation}
	Using Gaussian tail bounds and $\s>\sqrt{2}$, the latter integral in \eqref{B2.018} is bounded from above by
	\begin{equation}\label{B2.019}
		\frac{\sqrt{t}}{m(t)-\s t-B}\eee^{-(m(t)-\s t-B)^2/2t}.
	\end{equation}
	Inserting \eqref{B2.019} back into \eqref{B2.018}, there exists a constant $C>0$ such that \eqref{B2.018} is bounded from above by
	\begin{equation}\label{B2.0191}
		C t^{-1/2} \eee^{t-\s m(t)+\sigma^2 t/2- (m(t)-\s t-B)^2/2t}
		\leq C \eee^{(\sqrt 2 -\s)B}(1+o(1)).
	\end{equation}
	As $\sqrt 2-\s<0$, \eqref{B2.0191} converges to zero when first $t$ tends to infinity and then $B$ tends to infinity. This concludes the proof of Proposition \ref{Prop.B2.1}  in Case (a).
	
	Next, we prove Proposition \ref{Prop.B2.1} in Case (b).
	In this case, we bound the probability in \eqref{B2.016} from above by
	\be\Eq(B2.1)
	\P\left(\left|\sum_{k=1}^{n(t)} \eee^{\s x_k(t)+i\t y_k(t)-\s m(t)}\1_{x_k(t)<m(t)-B}\right|>\d\right),
	\ee
	which itself, using Lemma \thv(Lem.B1.1), is bounded from above by
	\be\Eq(B2.2)
	\P\left(\left|\sum_{k=1}^{n(t)} \eee^{\s x_k(t)+i\t y_k(t)-\s m(t)}\1_{x_k(t)<m(t)-B, x_k\in \mathcal{U}_{t,\g}}\right|>\d\right)+\e,
	\ee
	for all $C$ therein large enough. Using Chebyshev's inequality, we bound the probability in \eqv(B2.2) from above by
	\be\Eq(B2.3)
	\d^{-2}\E\left(\sum_{k,k'=1}^{n(t)}\eee^{\s x_k(t)+i\t y_k(t)
		+\s x_{k'}(t)-i\t y_{k'}(t)- 2\s m(t)}  \1_{x_k(t),x_{k'}(t)<m(t)-B, x_k,x_{k'}\in \mathcal{U}_{t,\g}}\right).
	\ee
	Using \eqref{cor.1}, the expectation in \eqref{B2.3} is equal to
	\begin{equation}\label{eq:extra.1}
		\eee^{- 2\s m(t)}\E\left(\sum_{k,k'=1}^{n(t)}\eee^{\l x_k(t)
			+\overline{\l} x_{k'}(t)+i\t \sqrt{1-\rho^2} (z_k(t)-z_{k'}(t))} \1_{x_k(t),x_{k'}(t)<m(t)-B, x_k,x_{k'}\in \mathcal{U}_{t,\g}}\right).
	\end{equation}
	As $Z$ is a variable speed branching Brownian motion on the same Galton-Watson tree as $X$, the expectation is equal to
	\begin{eqnarray}\Eq(B2.4)
		\E\bigg(\sum_{k,k'=1}^{n(t)}\eee^{\l x_k(t)
			+\overline{\l} x_{k'}(t)}\1_{x_k(t),x_{k'}(t)<m(t)-B, x_k,x_{k'}\in \mathcal{U}_{t,\g}} \E\big[\eee^{i\t \sqrt{1-\rho^2} (z_k(t)-z_{k'}(t))}\big]\bigg)\nonumber
		\\
		= \E\bigg(\sum_{k,k'=1}^{n(t)}\eee^{\l x_k(t)
			+\overline{\l} x_{k'}(t)}\1_{x_k(t),x_{k'}(t)<m(t)-B, x_k,x_{k'}\in \mathcal{U}_{t,\g}} 
		\eee^{\t^2(1-\rho^2)\left(t-tA\left(\frac{d(x_{k},x_{k'})}{t}\right)\right)}
		\bigg).
	\end{eqnarray}
	The expectation in \eqv(eq:extra.1) is bounded from above by
	\bea\Eq(B2.3.1)
	(I)+(II)&\equiv&\eee^{-2\s m(t)}\E\Bigg(\sum_{k,k'=1}^{n(t)}\eee^{\t^2(1-\rho^2)\left(t-tA\left(\frac{d(x_{k},x_{k'})} {t}\right)\right)}\eee^{\l x_k(t)+\overline{\l}x_{k'}(t)}
	\\ \nonumber
	&&\times\bigg( \1_{d(x_k,x_{k'})\leq r}
	+
	\1_{x_k(t),x_{k'}(t)<m(t)-B, x_k,x_{k'}\in \mathcal{U}_{t,\g},d(x_k,x_{k'})> r}\bigg)\Bigg).
	\eea
	(I) is equal to 
	\begin{equation}
		\eee^{-2\s m(t)}\int_0^{r}\dd s \eee^{\t^2(1-\rho^2)\left(t-tA(s/t)\right)}
		\eee^{2t-s}\int_{-\infty}^\infty \eee^{2\s x-\frac{x^2}{2tA(s/t)}} \left|\int_{-\infty}^{\infty} \dd y \eee^{-\l y-\frac{y^2}{2(t-tA(s/t))}}\right|^2.
	\end{equation}
	An explicit computation of the Gaussian integrals yields,
	\begin{multline}
		\eee^{-2\s m(t)}\int_0^{r}\dd s \eee^{\t^2(1-\rho^2)\left(t-tA(s/t)\right)} \eee^{2t-s} \eee^{2\sigma^2 tA(s/t)}\eee^{(\sigma^2-\rho^2\t^2)(t-tA(s/t))}\\
		\leq t^{3\s/\sqrt{2}}\eee^{((\sqrt{2}-\s)^2-\t^2)t}\int_0^{r}\dd s \eee^{-s+\sigma^2 tA(s/t)}
		\leq r \eee^{\sigma^2 r}t^{3\s/\sqrt{2}}\eee^{((\sqrt{2}-\s)^2-\t^2)t},    
	\end{multline}
	which converges to zero as $t\to\infty$ using that $\b\in B_2$.
	We bound (II) in \eqv(B2.3.1) from above by only keeping the localization constraint at the splitting time and further decomposing the event in the indicator function. It is bounded from above by
	\begin{equation}
		\E\left(\sum_{k,k'=1}^{n(t)}\eee^{\t^2(1-\rho^2)\left(t-tA\left(\frac{d(x_{k},x_{k'})}{t}\right)\right)}\eee^{\l x_k(t)+\overline{\l} x_{k'}(t)-2\s m(t)}\left( \1_{E_1}+\1_{E_2}+\1_{E_3}\right)\right)
	\end{equation}
	where
	\bea\Eq(B2.5)
	E_1&=&\left\{ r\leq d(x_k,x_{k'})\leq \left(\Sigma^2\right)^{-1}(t-B/2),
	x_k\left(d(x_k,x_{k'})\right)\leq  U_{A,\g}(d(x_k,x_{k'}))\right\},\nonumber\\
	E_2&=&\left\{d(x_k,x_{k'})> (\Sigma^2)^{-1}(t-B/2), \;x_k(d(x_k,x_{k'}))< U_{A,\g}(d(x_k,x_{k'}))-B/2\right\},\nonumber\\
	E_3&=&\left\{ \Bracket{d(x_k,x_{k'})> (\Sigma^2)^{-1}(t-B/2), \;x_k(t),x_{k'}(t)<m(t)-B,}{x_k(d(x_k,x_{k'}))> U_{A,\g}(d(x_k,x_{k'}))-B/2}\right\}.
	\eea
	Using the many-to-two lemma,
	\begin{equation}
		\E\left(\sum_{k,k'=1}^{n(t)}\eee^{\t^2(1-\rho^2)\left(t-tA\left(\frac{d(x_{k},x_{k'})}{t}\right)\right)}\eee^{\l x_k(t)+\overline{\l} x_{k'}(t)-2\s m(t)}\1_{E_1}\right)
	\end{equation}
	is equal to
	\bea\Eq(B2.6)
	&&\int_{r}^{\left(\Sigma^2\right)^{-1}(t-B/2)} \dd s \eee^{2t-s-2\s m(t)} \int_{-\infty}^{U_{A,\g}(s)}\frac{\dd x}{\sqrt{2\pi tA(s/t)}} \eee^{2\s x-\frac{x^2}{2 tA(s/t)}}\\\nonumber
	&&\times \int_{-\infty}^{\infty}\frac{\dd y \eee^{\l y-\frac{y^2}{2t\left(1-A(s/t)\right)}}}{\sqrt{2\pi t\left(1-A(s/t)\right)}} 
	\int_{-\infty}^{\infty}\frac{\dd y' \eee^{\overline{\l} y'-\frac{(y')^2}{2t\left(1-A(s/t)\right)}}}{\sqrt{2\pi t\left(1-A(s/t)\right)}}
	.
	\eea
	Proceeding as in \eqv(B1.14)-\eqv(B1.17), we bound \eqv(B2.6) from above by
	\bea\Eq(B2.7)
	&&\eee^{2t+\left(\s^2-\t^2\right)t-2\s m(t)}\bigg(\int_{\left(\Sigma^2\right)^{-1}(r)}^{t-B/2} \dd s \eee^{-s+(\s^2+\t^2)tA(s/t)}\1_{U_{A,\g}(s)-2\s t A(s/t)\geq 0}\\
	&&+ \int_{r}^{\left(\Sigma^2\right)^{-1}(t-B/2)} \dd s \eee^{-s+(\s^2+\t^2)tA(s/t)}\1_{U_{A,\g}(s)-2\s t A(s/t)< 0} \frac{\eee^{-\frac{\left(U_{A,\g}(s)-2\s t A(s/t)\right)^2}{2tA(s/t)}}}{\sqrt{2s}-2\s \sqrt{tA(s/t)}} \bigg).\nonumber
	\eea
	The integrand in the second summand in \eqv(B2.7) is
	\bea\Eq(B2.8)
	&&\frac{t^{\frac{\s}{\sqrt{2}}}}{\sqrt{2s}-2\s \sqrt{tA(s/t)}}\eee^{(\sqrt{2}-\s)^2t-\t^2t-\left(\sqrt{2s}-\s\sqrt{tA(s/t)}\right)^2+\t^2 tA(s/t)}
	\\ \Eq(B2.8.1)
	&&\qquad\qquad \times \eee^{\left(2\s-\sqrt{\frac{2s}{tA(s/t)}}\right)\left(\left(tA(s/t)\land (t-tA(s/t))\land C\right)^\g-\frac{\sqrt{tA(s/t)}}{2\sqrt{2s}}\log(t)\right)}.
	\eea
	
	We rewrite the exponential in \eqref{B2.8} as
	\begin{equation}\label{B2.8.2}
		\eee^{\left((\sqrt{2}-\s)^2-\t^2\right)(t-tA(s/t)) }.
	\end{equation}
	
	Hence, for $\g<\alpha< 1$, we can bound the second integral in \eqref{B2.7} from above by
	\be\Eq(B2.9)
	\int_{r}^{\left(\Sigma^2\right)^{-1}(t-t^\a)} \dd s \frac{t^{\frac{\s}{\sqrt{2}}}}{\sqrt{2s}-2\s \sqrt{tA(s/t)}}\eee^{\left((\sqrt{2}-\s)^2-\t^2\right)(t-tA(s/t)) +O(s^\g) }\to 0
	\ee
	as $t\uparrow \infty$. Moreover, the integrand of the first integral in \eqref{B2.7} is bounded from above by
	\bea\Eq(B2.10)
	\left(1+o(1)\right) 
	\eee^{\left((\sqrt{2}-\s)^2-\t^2\right)(t-tA(s/t)) }\eee^{\left(\sqrt{\frac{2s}{tA(s/t)}}-2\s\right)\left((t-tA(s/t))\vee C\right)^\g },
	\eea
	as the polynomial terms in $t$ cancel exactly. 
	As $(t-tA(s/t))\vee C= (t-tA(s/t))$ for all such s, for B sufficiently large, we get
	\begin{equation}\label{B2.101}
		\int_{\left(\Sigma^2\right)^{-1}(t-t^\a)}^{\left(\Sigma^2\right)^{-1}(t-B/2)}\dd s \eee^{\left((\sqrt{2}-\s)^2-\t^2\right)(t-tA(s/t)) }\eee^{\left(\sqrt{\frac{2s}{tA(s/t)}}-2\s\right)\left((t-tA(s/t))\right)^\g }\leq \eee^{-cB},
	\end{equation}
	for some constant $c>0$ for all $B$ large enough. 
	Using the bounds \eqref{B2.9} and \eqref{B2.101} in \eqref{B2.7}, we bound the first summand in \eqv(B2.7) from above by
	\be\Eq(B2.11)
	t^{\frac{\s}{\sqrt{2}}}\eee^{\left(\sqrt{2}-\s\right)^2t-\t^2 t}\int_{r}^{t-B/2} \dd s \eee^{-s+(\s^2+\t^2)tA(s/t)}\1_{U_{A,\g}(s)-2\s t A(s/t)\geq 0}.
	\ee
	Observe that 
	\begin{equation}\label{B2.11.2}
		\eee^{-s+(\s^2+\t^2)tA(s/t)}\1_{U_{A,\g}(s)-2\s t A(s/t)\geq 0}\leq
		\eee^{\left(\sqrt{2}-\s\right)^2t-\t^2 t+ O\left(\left(((tA(s/t)\land (t-tA(s/t)))\vee C\right)^\g\right) }.
	\end{equation}
	Using \eqref{B2.11.2} and that $\t^2tA(s/t))\leq s/2$ as $\t>1/\sqrt{2}$ in Case (b),  the integral in \eqv(B2.11) is bounded by a constant. Hence, \eqref{B2.11} converges to zero as $t$ tends to infinity.
	Next, we bound 
	\begin{equation}
		\E\left(\sum_{k,k'=1}^{n(t)}\eee^{\t^2(1-\rho^2)\left(t-tA\left(\frac{d(x_{k},x_{k'})}{t}\right)\right)}\eee^{\l x_k(t)+\overline{\l} x_{k'}(t)-2\s m(t)}\1_{E_2}\right)
	\end{equation}
	from above by
	\bea\Eq(B2.12)
	&&\int_{\left(\Sigma^2\right)^{-1}(t-B/2)}^t \dd s \eee^{2t-s-2\s m(t)} \int_{-\infty}^{U_{A,\g}(s)-B/2}\frac{\dd x}{\sqrt{2\pi tA(s/t)}} \eee^{2\s x-\frac{x^2}{2 tA(s/t)}}\\\nonumber
	&&\times \int_{-\infty}^{\infty}\frac{\dd y}{\sqrt{2\pi t\left(1-A(s/t)\right)}} \eee^{\b y-\frac{y^2}{2t\left(1-A(s/t)\right)}}
	\int_{-\infty}^{\infty}\frac{\dd y'}{\sqrt{2\pi t\left(1-A(s/t)\right)}} \eee^{\overline{\b} y'-\frac{(y')^2}{2t\left(1-A(s/t)\right)}}.
	\eea
	Proceeding as in the second moment computation above, and using that  $U_{A,\g}(s)-B/2-2\s t A(s/t)< 0$ for all $t$ large enough on the domain of integration, we can bound  \eqv(B2.12) from above by
	\be\Eq(B2.13)
	\eee^{2t+\left(\s^2-\t^2\right)t-2\s m(t)} \int_{\left(\Sigma^2\right)^{-1}(t-B/2)}^t \dd s \eee^{-s+(\s^2+\t^2)tA(s/t)} \frac{\eee^{-\frac{\left(U_{A,\g}(s)-B/2-2\s t A(s/t)\right)^2}{2tA(s/t)}}}{\sqrt{2s}-2\s \sqrt{tA(s/t)}} .
	\ee
	Finally, we bound 
	\begin{equation}
		\E\left(\sum_{k,k'=1}^{n(t)}\eee^{\t^2(1-\rho^2)\left(t-tA\left(\frac{d(x_{k},x_{k'})}{t}\right)\right)}\eee^{\l x_k(t)+\overline{\l} x_{k'}(t)-2\s m(t)}\1_{E_3}\right)
	\end{equation}
	from above by
	\bea\Eq(B2.015)
	&&\int_{\left(\Sigma^2\right)^{-1}(t-B/2)}^t \dd s \eee^{2t-s-2\s m(t)} \int_{U_{A,\g}(s)-B/2}^{U_{A,\g}}\frac{dx}{\sqrt{2\pi tA(s/t)}} \eee^{2\s x-\frac{x^2}{2 tA(s/t)}} \nonumber\\
	&&\times \int_{-\infty}^{m(t)-B-x}\frac{\dd y \eee^{\b y-\frac{y^2}{2t\left(1-A(s/t)\right)}}}{\sqrt{2\pi t\left(1-A(s/t)\right)}} 
	\int_{-\infty}^{m(t)-B-x}\frac{\dd y' \eee^{\overline{\b} y'-\frac{(y')^2}{2t\left(1-A(s/t)\right)}}}{\sqrt{2\pi t\left(1-A(s/t)\right)}} 
	.
	\eea
	As $A(s/t) \leq 1 $, it holds on the domain of integration of $x$ that
	\be\Eq(B2.15)
	m(t)-B-x
	\leq m(t)-U_{A,\g}(s)-B/2 
	\leq\sqrt{2}(t-s)-B/2.
	\ee
	Hence, after also completing the squares, the integrals with respect to $y$ and $y'$ are bounded from above by
	\bea\Eq(B2.14)
	&& \eee^{(\s^2-\t^2)t\left(1-A(s/t)\right)}\int_{-\infty}^{\sqrt{2}(t-s)-B/2}\frac{\dd y}{\sqrt{2\pi t\left(1-A(s/t)\right)}} \eee^{ -\frac{\left(y-\b t(1-A(s/t))\right)^2}{2t\left(1-A(s/t)\right)}}
	\nonumber\\
	&&\times\int_{-\infty}^{\sqrt{2}(t-s)-B/2}\frac{\dd y'}{\sqrt{2\pi t\left(1-A(s/t)\right)}} \eee^{ -\frac{\left(y'-\overline\b t(1-A(s/t))\right)^2}{2t\left(1-A(s/t)\right)}}
	\nonumber\\
	&\leq& \eee^{(\s^2-\t^2)t\left(1-A(s/t)\right)} \eee^{-2\frac{\left(\left(\sqrt{2}-\s\right)B/2 -B/2\right)^2}{2t(1-A(s/t))}}\leq \eee^{(\s^2-\t^2)t\left(1-A(s/t)\right)} \eee^{-\frac{\left(\left(\sqrt{2}-\s\right) -1\right)^2B}{2}}.
	\eea
	The $x$ integral can now be bounded as before, but we get a negative exponential in $B$, which gives us an extra decay, which concludes the proof.
\end{proof}
\begin{proof}[Proof of Theorem~\ref{phaseB2.1}:]
	Denote by $\mathbb{M}$ the space of locally finite counting measures on $\mathbb{R}\cup \{+\infty\}$ endowed with the vague topology. Consider for $B>0$ the functional $\Phi_{\beta, B}\colon \mathbb{M}\rightarrow \mathbb{R}$ defined by
	\begin{align}
		\Phi_{\beta, B}\left(\sum_{i\in I}\delta_{x_i}\right)=\sum_{i\in I} \eee^{\beta x_i}\mathbbm{1}_{x_i>-B}
	\end{align}
	where $I$ is a countable index set and $\sum_{i\in I}\delta_{x_i}$ a locally finite counting measure. The set of locally finite counting measures on which the functional $\Phi_{\beta,B}$ is not continuous (the measure charging $-B$ or $+\infty$) has measure zero w.r.t.~to the law of $\mathcal{E}_{\sigma_b,\sigma_e}$. Hence, by the continuous mapping theorem, it follows that in law and as $t\uparrow \infty$,
	\begin{align}
		\Phi_{\beta, B}\left(\mathcal{E}_{\sigma_b,\sigma_e}(t)\right)\Rightarrow \Phi_{\beta, B}\left(\mathcal{E}_{\sigma_b,\sigma_e}\right)=\sum_{k,l}\eee^{\beta (\eta_k +\Delta^{(k)}_l)} \mathbbm{1}_{\eta_k+\Delta^{(k)}_l>-B}.
	\end{align}
	Note that by Proposition \ref{Prop.B2.1}, for all $\epsilon>0$ and $\delta>0$, there exists $B_0$ such that for all $B\geq B_0$ and all $t$ large enough,
	\begin{align}
		\P\left(|\XX_{\beta,1}(t)-\Phi_{\beta,B}(\mathcal{E}_{\sigma_b, \sigma_e}(t)|>\delta \right) <\epsilon.
	\end{align}
	Hence, by Slutsky's theorem, $\XX_{\beta,1}(t)$ converges in law to
	\begin{align}        \XX_{\beta,1}=\lim_{B\uparrow\infty}\sum_{k,l}\eee^{\beta(\eta_k+\Delta^{(k)}_l)}\mathbbm{1}_{\eta_k+\Delta^{(k)}_l>-B}.
	\end{align}
	
\end{proof}

\begin{proof}[Proof of Theorem~\ref{phaseB2.2}:]
	Using representation \eqref{cor.1}, $e^{-\s m(t)}\XX_{\beta,\rho}(t)$ is equal in law to
	\begin{align}
		\sum_{k\leq n(t)}\eee^{(\s+i\rho\t)(x_k(t)-m(t))+i\tau \sqrt{1-\rho^2}z_k(t)+i\rho\t m(t)},
	\end{align}
	where $\left\{z_k(t), \, k=1,\dotsc, n(t)\right\}$ is an independent copy of the variable-speed BBM $X(t)$ defined on the identical Galton-Watson tree. If $\rho\notin \{-1,1\}$, then by \cite[Lemma~3.2]{MRV13} (using that every path in variable-speed BBM can be seen as a time-changed path of standard BBM) and an adaption of the following discussion after the statement therein, which handles standard BBM, to variable-speed BBM, we get that
	\begin{align}
		\mathcal{G}(t):=\sum_{k\leq n(t)} \delta_{(x_k(t)-m(t),\exp[i\sqrt{1-\rho^2}\t z_k(t)+i\rho \t m(t)])}
	\end{align}
	converges weakly as $t\uparrow \infty$ to
	\begin{align}
		\mathcal{G}:=\sum_{k,l\geq 1}\delta_{(\eta_k+\Delta^{(k)}_l, U^{(k)}\tilde{W}_l^{(k)})},
	\end{align}
	where $\eta_k, \Delta^{(k)}_l$ are as in \eqref{eq:extremal_process}, $(U^{(k)})_{k\geq 1}$ are i.i.d.\ uniformly distributed on the unit circle and $\tilde{W}^{(k)}_l$ the atoms of a point process on the unit circle. We denote by $\tilde{\mathbb{M}}$ the space of locally finite counting measures on $\{\mathbb{R}\cup \{+\infty\}\}\times \{ z\in \mathbb{C}:|z|=1\}$, endowed with the topology of vague convergence. For $B>0$, consider the functional $\tilde{\Phi}_{\beta,B}\colon\tilde{\mathbb{M}}\rightarrow \mathbb{C}$ defined as
	\begin{align}
		\tilde{\Phi}_{\beta,B}\left(\sum_{k\in I}\delta_{(x_k,z_k)}\right):=\sum_{k\in I}\eee^{\beta x_k}z_k \mathbbm{1}_{\eta_k +\Delta^{(k)}_l\geq-B},
	\end{align}
	where $I$ is a countable index set. The set of locally finite counting measures on $\tilde{\mathbb{M}}$ on which the functional $ \tilde{\Phi}_{\beta,B}$ is not continuous, i.e., that charge $(-B,\cdot)$ or $(+\infty,\cdot)$, has measure zero w.r.t.~the law of $\mathcal{G}$. Hence, by the continuous mapping theorem, it follows that, in law,
	\begin{align}\label{eq:phaseB2.2.1}
		\lim_{t\uparrow \infty}\tilde{\Phi}_{\beta,B}(\mathcal{G}_t)= \tilde{\Phi}_{\beta,B}(\mathcal{G})=\sum_{k,l\geq 1} \eee^{(\sigma+i\rho\t)(\eta_k+\Delta^{(k)}_l)}U^{(k)}\tilde{W}^{(k)}_l\mathbbm{1}_{\eta_k+\Delta^{(k)}_l\geq -B}.
	\end{align}
	Now, $e^{i\rho\t(\eta_k+\Delta^{(k)}_l)}U^{(k)}$ is uniformly distributed on the unit circle and so \eqref{eq:phaseB2.2.1} is equal in law to
	\begin{align}
		\sum_{k,l\geq 1}\eee^{\sigma(\eta_k+\Delta^{(k)}_l)}U^{(k)}\tilde{W}^{(k)}_l\mathbbm{1}_{\eta_k+\Delta^{(k)}_l\geq -B}.
	\end{align}
	By Proposition~\ref{Prop.B2.1}, for all $\epsilon,\delta>0$ there exists $B_0$ such that for all $B\geq B_0$ and $t$ large,
	\begin{align}
		\P\left(\bigg| \XX_{\beta,\rho}(t)-\tilde{\Phi}_{\s+i\rho \t,B}(\mathcal{G}(t) \bigg|>\delta\right)<\epsilon.
	\end{align}
	Hence, by Slutsky's theorem, $\XX_{\beta,\rho}(t)$ converges in law, as $t\uparrow\infty$, to
	\begin{align}
		\XX_{\beta,\rho}=\lim_{B\uparrow\infty}\sum_{k,l\geq 1}\eee^{\sigma(\eta_k+\Delta^{(k)}_l)}U^{(k)}\tilde{W}^{(k)}_l\mathbbm{1}_{\eta_k+\Delta^{(k)}_l\geq -B}=\sum_{k,l\geq 1}\eee^{\sigma(\eta_k+\Delta^{(k)}_l)}U^{(k)}\tilde{W}^{(k)}_l.
	\end{align}
	The latter we may rewrite as
	\begin{align}
		\sum_{k\geq 1}e^{\sigma\eta_k}U^{(k)}{W}^{(k)},
	\end{align}
	where ${W}^{(k)}:=\sum_{l\geq 1}e^{\sigma\Delta^{(k)}_l}\tilde{W}^{(k)}_l$ are i.i.d.~random variables. Thus, conditionally on $Y_{\s_b}$, the law of $\XX_{\beta,\rho}$ is complex isotropic $(\sqrt{2}/\s)$-stable.
\end{proof}

\section{Phase B3}
Let
\begin{equation}
	N_{\t,\s}(t)= \eee^{-t(1/2 +\s^2)} \XX_{\b,\rho}.
\end{equation}
The proof of Theorem \ref{clt.B3} relies on the method of moments similar to the proof of Theorem 1.11 in \cite{hartklim18}. In a first step, we prove in Lemma \ref{LemB3.1} that the limit of the second moment as $t\to\infty$ is finite. Then, we show in Lemmas \ref{LemB3.2} and \ref{LemB3.3} that we can introduce a suitable truncation. Finally, in Lemma \ref{Lem.const3} we study the asymptotics of all moments of the truncated version of $N_{\t,\s}(t)$.
\begin{lemma}[The second moment]\label{LemB3.1} 
	For $\b\in B_3$, there exists a constant $C$ such that
	\begin{equation}\label{B3.1}
		\lim_{t\to \infty}\E(| N_{\t,\s}(t)|^2)=C.
	\end{equation}
\end{lemma}
\begin{proof}
	By \eqref{cor.1} the expectation in \eqref{B3.1} is equal to 
	\begin{equation} \label{B3.2} \E\left(\sum_{k,k'=1}^{n(t)}\eee^{\s x_k(t)+i\t y_k(t)
			+\s x_{k'}(t)-i\t y_{k'}(t)- t(1+2\s^2)} \right).
	\end{equation}
	Using that the correlation coefficient of $X$ and $Y$ is $\rho$ and setting $\lambda=(\s+i\t\rho)$, \eqref{B3.2} can be written as
	\begin{eqnarray} \label{B3.2.1} &&\E\left(\sum_{k,k'=1}^{n(t)}\eee^{\l x_k(t)
			+\overline{\l}  x_{k'}(t)+i\sqrt{1-\rho^2}\t (z_k(t) -z_{k'}(t))- t(1+2\s^2)} \right)\nonumber\\
		&&=\E\left(\sum_{k,k'=1}^{n(t)} \eee^{-t(1+2\s^2)-(t-tA(d(x_k,x_{k'})/t))(1-\rho^2) \t^2}\eee^{\l x_k(t)
			+\overline{\l}  x_{k'}(t)}\right)
		,
	\end{eqnarray}
	where we used the independence of $X$ and $Z$ when conditioning on the underlying tree. By the many-to-two lemma, the expectation in \eqref{B3.2.1} is equal to
	\begin{eqnarray}\label{B3.3}
		&& K  \eee^{-t(1+2\s^2)}\int_0^{t}\dd s \eee^{-\t^2(1-\rho^2)\left(t-tA(s/t)\right)}
		\eee^{2t-s}\int_{-\infty}^\infty \frac{\dd x}{\sqrt{2\pi tA(s/t)}} \eee^{2\s x-\frac{x^2}{2tA(s/t)}} \nonumber\\
		&& \times\left|\int_{-\infty}^{\infty} \frac{
			\dd y}{\sqrt{2\pi(t-tA(s/t))}} \eee^{\l y-\frac{y^2}{2(t-tA(s/t))}}\right|^2\nonumber\\
		&&= K\eee^{-t(1+2\s^2)}\int_0^{t}\dd s \eee^{2t-s+(\s^2-\t^2)\left(t-tA(s/t)\right)}
		\int_{-\infty}^\infty \frac{\dd x}{\sqrt{2\pi tA(s/t)}} \eee^{2\s x-\frac{x^2}{2tA(s/t)}} .
	\end{eqnarray}
	As the integral with respect to $x$ evaluates to $\eee^{2\s^2 tA(s/t)}$, we get in total
	\begin{equation}\label{B3.4}
		K  \int_0^t\dd s \eee^{t-s - (\s^2 +\t^s)(t-tA(s/t))}. 
	\end{equation}
	Now, for $s$ close to $t$, we have
	\begin{equation}\label{B3.5}
		t-tA(s/t)= \s_e^2 (t-s)+ O(s^2/t).
	\end{equation}
	We note that $tA(s/t)<s$ and $\s^2+\t^2<1$. Using this and \eqref{B3.5} we deduce that the integral in \eqref{B3.4} is asymptotically equal to
	\begin{eqnarray}\label{B3.6}
		\lim_{R\to \infty}  K  \int_{t-R}^t \dd s \eee^{t-s - (\s^2 +\t^s)(t-tA(s/t))}= \lim_{R\to \infty}  K  \int_{t-R}^t \dd s \eee^{t-s - (\s^2 +\t^s) \s_e^2 (t-s)}\nonumber\\
		=\frac{K}{(\s^2 +\t^2) \s_e^2-1}.
	\end{eqnarray}
\end{proof}

\begin{lemma}\label{LemB3.2}
	Let $\b \in B_3$. For all $\e>0$ and $\d>0$, there exists $B_0$ such that for all $B>B_0$, such that for all $t$ large enough,
	\begin{equation}\label{B3.7}
		\P\left(\left|
		\sum_{k=1}^{n(t)} \eee^{\s x_k(t)+i\t y_k(t)-\left(\frac{1}{2}+\s^2\right)t}\1_{x_k(t)>2\s t+B\sqrt{t}}
		\right|>\d \right)<\e.
	\end{equation}
\end{lemma}
\begin{proof}
	Using Chebyshev's inequality, the probability in \eqref{B3.7} is smaller than
	\begin{equation}\label{B3.8}
		\d^{-2}\E\left(
		\left|
		\sum_{k=1}^{n(t)} \eee^{\s x_k(t)+i\t y_k(t)-\left(\frac{1}{2}+\s^2\right)t}\1_{x_k(t)>2\s t+B\sqrt{t}}
		\right|^2\right).
	\end{equation}
	Computing the square and using as in \eqref{B3.2.1} that the correlation coefficient between $X$ and $Y$ is equal to $\rho$ and $\lambda=(\s+i\t \rho)$, the expectation in \eqref{B3.8} is equal to
	\begin{equation}\label{B3.9}
		\E\left(\sum_{k,k'=1}^{n(t)} \eee^{-t(1+2\s^2)-(t-tA(d(x_k,x_{k'})/t))(1-\rho^2) \t^2}\eee^{\l x_k(t)
			+\overline{\l}  x_{k'}(t)}\1_{x_k(t),x_{k'}(t)>2\s t+B\sqrt{t}}\right) .
	\end{equation}
	We split the expectation in \eqref{B3.9} into two summands, (I) and (II), distinguishing whether $d(x_k,x_{k'})<t-r$ or not (we choose $r>0$ later in the proof).
	\begin{equation}\label{B3.10}
		(I)\leq \E\left(\sum_{k,k'=1}^{n(t)}\1_{d(x_k,x_{k'})\leq t-r} \eee^{-t(1+2\s^2)-(t-tA(d(x_k,x_{k'})/t))(1-\rho^2) \t^2}\eee^{\l x_k(t)
			+\overline{\l}  x_{k'}(t)}\right),
	\end{equation}
	which is equal to \eqref{B3.3} when only integrating $s$ from $0$ to $t-r$. Following the computations up to \eqref{B3.6}, this can be made smaller than $\e/3$ by choosing $r$ sufficiently large.
	(II) is equal to
	\begin{eqnarray}\label{B3.11}
		&&   K  \eee^{-t(1+2\s^2)}\int_{t-r}^{t}\dd s \eee^{-\t^2(1-\rho^2)\left(t-tA(s/t)\right)}
		\eee^{2t-s}\int_{-\infty}^\infty \frac{\dd x}{\sqrt{2\pi tA(s/t)}} \eee^{2\s x-\frac{x^2}{2tA(s/t)}} \nonumber\\
		&&\times\left| \int_{2\s t+B\sqrt{t}-x}^{\infty} \frac{
			\dd y}{\sqrt{2\pi(t-tA(s/t))}} \eee^{\l y-\frac{y^2}{2(t-tA(s/t))}}\right|^2.
	\end{eqnarray}
	For any $R>0$, if $2\s t +B\sqrt{t}-x>R$, then the absolute value of the integral with respect to $y$, can be made arbitrarily small (for any finite $r$) by choosing $R$ large enough (as $t-tA(s/t)$ is by \eqref{B3.5} of order $r$). We fix such $R>0$. Finally,
	we look at the complement $2\s t +B\sqrt{t}-x<R$. On this event, \eqref{B3.11} is smaller than
	\begin{eqnarray}\label{B3.12}
		&&   K  \eee^{-t(1+2\s^2)}\int_{t-r}^{t}\dd s \eee^{-\t^2(1-\rho^2)\left(t-tA(s/t)\right)}
		\eee^{2t-s}\int_{2\s t + B\sqrt{t}-R}^\infty \frac{\dd x}{\sqrt{2\pi tA(s/t)}} \eee^{2\s x-\frac{x^2}{2tA(s/t)}} \nonumber\\
		&&\times\left| \int_{-\infty}^{\infty} \frac{
			\dd y}{\sqrt{2\pi(t-tA(s/t))}} \eee^{\l y-\frac{y^2}{2(t-tA(s/t))}}\right|^2.
	\end{eqnarray}
	We compute the integral with respect to $y$ as before and get that \eqref{B3.12} is equal to
	\begin{equation}
		K \eee^{-t(1+2\s^2)}\int_{t-r}^{t}\dd s \eee^{(\s^2-\t^2)\left(t-tA(s/t)\right)}
		\eee^{2t-s}\int_{2\s t +B\sqrt{t}-R}^\infty \frac{\dd x \eee^{-\frac{(x-2\s tA(s/t))^2}{2tA(s/t)}}\eee^{2\s^2 tA(s/t)}}{\sqrt{2\pi tA(s/t)}}.
	\end{equation}
	We observe that the integral
	\begin{equation}\label{B3.13}
		\int_{2\s t +B\sqrt{t}-R}^\infty \frac{\dd x \eee^{-\frac{(x-2\s tA(s/t))^2}{2tA(s/t)}}}{\sqrt{2\pi tA(s/t)}} 
	\end{equation}
	can be made arbitrarily small by choosing $B$ large enough, as $tA(s/t)=t-O(r)$ (see \eqref{B3.5}).
	This implies, by the second moment computation in Lemma \ref{LemB3.1} that \eqref{B3.12} can be made as small as we like by taking $B$ large enough. This together with our previous bounds concludes the proof of Lemma \ref{LemB3.2}.
\end{proof}
Set
\be\Eq(B3.22)
\mathcal{S}_{r,t,\g} : =\left\{x(\cdot)\in C\left([0,t],\R\right) :\forall s\in [(\Sigma^2)^{-1}(r),t]: x_k(s)> 2\s tA(s/t)+(tA(s/t))^\g
\right\}
.
\ee
\begin{lemma}\label{LemB3.3}
	Let $\b \in B_3$, $\g>1/2$ and $B>0$.  For all $\e>0$ and $\d>0$, there exists $r_0$ such that for all $r>r_0$ and uniformly for all $t$ large enough,
	\begin{equation}\label{B3.14}
		\P\left(\left|
		\sum_{k=1}^{n(t)} \eee^{\s x_k(t)+i\t y_k(t)-\left(\frac{1}{2}+\s^2\right)t}\1_{x_k(t)<2\s t+B\sqrt{t}, x_k\not\in \mathcal{S}_{r,t,\g}}
		\right|>\d \right)<\e.
	\end{equation}  
\end{lemma}
\begin{proof}
	As in the proof of Lemma \ref{LemB3.2}, we use Chebyshev's inequality to bound the probability in \eqref{B3.14} from above by 
	\begin{equation}\label{B3.15}  
		\d^{-2}\E\Big(
		\Big|
		\sum_{k=1}^{n(t)} \eee^{\s x_k(t)+i\t y_k(t)-\left(\frac{1}{2}+\s^2\right)t}\1_{x_k(t)\leq 2\s t+B\sqrt{t},x_k\not\in \mathcal{S}_{r,t,\g}}
		\Big|^2\Big).
	\end{equation}
	Using the same reformulation as in \eqref{B3.9} and dropping part of the restrictions for one of the particles, the expectation in \eqref{B3.15} is bounded from above by $\eee^{-t(1+2\s^2)}$ times
	\begin{equation}\label{B3.16}
		\E\Big(\sum_{k,k'=1}^{n(t)} \eee^{-(t-tA(d(x_k,x_{k'})/t))(1-\rho^2) \t^2}\eee^{\l x_k(t)
			+\overline{\l}  x_{k'}(t)}\1_{x_k(t),x_{k'}(t)\leq 2\s t+B\sqrt{t}, x_k\not\in \mathcal{S}_{r,t,\g}}\Big).
	\end{equation}
	We have seen in the proof of Lemma \ref{B3.1} that
	\begin{equation}\label{B3.17}
		\E\Big(\sum_{k,k'=1}^{n(t)} \eee^{-t(1+2\s^2)-(t-tA(d(x_k,x_{k'})/t))(1-\rho^2) \t^2}\eee^{\l x_k(t)
			+\overline{\l}  x_{k'}(t)}\1_{d(x_k,x_{k'})< t-r}\Big)
	\end{equation}
	can be made arbitrarily small by choosing $r$ large enough. Moreover,
	\begin{eqnarray}\label{B3.18}
		&\E\Big(\sum_{k,k'=1}^{n(t)} \eee^{-t(1+2\s^2)-(t-tA(d(x_k,x_{k'})/t))(1-\rho^2) \t^2}\nonumber \\ &\qquad\times \eee^{\l x_k(t)
			+\overline{\l}  x_{k'}(t)}\1_{d(x_k,x_{k'})> t-r}\1_{x_k(t)-x_k(d(x_k,x_{k'}))<-R}\Big)
	\end{eqnarray}
	can be made arbitrarily small by choosing $R$ large enough. Hence, it suffices to bound 
	\begin{eqnarray}\label{B3.19}
		&&\E\Big(\sum_{k,k'=1}^{n(t)} \eee^{-t(1+2\s^2)-(t-tA(d(x_k,x_{k'})/t))(1-\rho^2) \t^2}\eee^{\l x_k(t)
			+\overline{\l}  x_{k'}(t)}\1_{x_k(t)\leq 2\s t+B\sqrt{t}}\\ \nonumber
		&&\qquad\times  \1_{x_k\not\in \mathcal{S}_{r,t,\g}}\1_{d(x_k,x_{k'})> t-r}\1_{x_k(t)-x_k(d(x_k,x_{k'}))>-R}\Big).
	\end{eqnarray}
	We observe that $x_k(t)-x_k(d(x_k,x_{k'}))>-R$ and $x_{t}<2\s t+B\sqrt{t}$ imply that 
	\begin{equation}\label{B3.20}
		x_k(d(x_k,x_{k'}))\leq 2\s t+B\sqrt{t}+R.
	\end{equation}
	We distinguish two cases: In the first case, we find $s\in [(\Sigma^2)^{-1}(r),d(x_k,x_{k'})]$ and in the second case, there exists $s\in [d(x_k,x_{k'}),t]$, such that the particle position at time $s$ is too high. We start with the first case. On this event, we can bound \eqref{B3.19} from above by $K  e^{-t(1+2\s^2)}$ times
	\begin{eqnarray}\label{B3.21}
		&&\int_{t-r}^{t}\dd q \eee^{-\t^2(1-\rho^2)\left(t-tA(q/t)\right)}
		\eee^{2t-q} \E\left( \eee^{\l y_1(t-tA(q/t)) +\overline{\l} y_2(t-tA(q/t)) }\right)\\ \nonumber
		&& \times\E\left(e^{2 \s x_1(tA(q/t))}\1_{\exists s\in [(\Sigma^2)^{-1}(r),q]: x_1(tA(s/t))\leq 2\s tA(s/t)+(tA(s/t))^\g , x_1(tA(q/t))<2\s t+B\sqrt{t}+R.}\right)
		,
	\end{eqnarray}
	where $x_1$ is a standard Brownian motion and $y_1(t-tA(q/t)),y_2(t-tA(q/t))$ are independent centered Gaussian random variables with variance $t-tA(q/t)$.
	Noting that 
	\begin{equation}\label{B3.222}
		\zeta(s)= x_1(s)-\frac{s}{tA(q/t) }x_1(tA(q/t)) ,
	\end{equation}
	is a standard Brownian bridge from zero to zero in time $tA(q/t)$, we see that the second expectation in \eqref{B3.21} is bounded from above by
	\begin{eqnarray}\label{B3.24}
		&& \E\left(\eee^{2 \s x_1(tA(q/t))}\1_{\exists s\in [r,tA(q/t)]: \zeta(s)\leq 2\s s+s^\g-\frac{s}{tA(q/t)}(2 \s t +B\sqrt{t}+R) }\right) \nonumber\\
		&&\leq 
		\E\left(\eee^{2 \s x_1(tA(q/t))}\1_{\exists s\in [r,tA(q/t)]: \zeta(s)\leq s^\g-O\left(\frac{s}{\sqrt{t}}\right)-O\left(\frac{(r\lor R)s}{t}\right) }\right) =
		\E\left(\eee^{2 \s x_1(tA(q/t))}\right)\nonumber \\ 
		&&\qquad\times \P\left(\exists s\in [r,tA(q/t)]: \zeta(s)\leq s^\g-O\left(\frac{s}{\sqrt{t}}\right)-O\left(\frac{(r\lor R)s}{t}\right) \right),
	\end{eqnarray}
	where in the last step we use that
	$\zeta$ is independent from $x_1(tA(q/t))$.  For $r$ large enough, by standard Brownian bridge estimates, the probability in \eqref{B3.24} smaller than $\e/C$, for some constant $C>1$ which is independent of $\e$ (and can be chosen arbitrarily large by increasing $r$).
	Plugging this back into \eqref{B3.21} we get that \eqref{B3.21} is bounded from above by 
	\begin{equation}\label{B3.25}
		\e/C\int_{t-r}^{t}\dd q \eee^{-\t^2(1-\rho^2)\left(t-tA(q/t)\right)}
		\eee^{2t-q} \E\left( \eee^{\l y_1(t-tA(q/t)) +\overline{\l} y_2(t-tA(q/t)) }\right)
		\E\left(e^{2 \s x_1(tA(q/t))}\right).
	\end{equation}
	In the second case,  where $s\in [d(x_k,x_{k'}),t]$, we bound \eqref{B3.19} from above by $Ke^{-t(1+2\s^2)}$ times 
	\begin{eqnarray}
		\label{B3.26}
		&&\int_{t-r}^{t}\dd q \eee^{-\t^2(1-\rho^2)\left(t-tA(q/t)\right)}
		\eee^{2t-q} \E\left(\eee^{\overline{\l} y_2(t-tA(q/t)) }\right)\E\left(\eee^{2 \s x_1(tA(q/t))}\right)\nonumber\\
		&&\quad\times
		\E\left( \eee^{\l (y_1(t)-y_1(tA(q/t)))}\1_{E_s}\right),
	\end{eqnarray}
	where $y_1$ is a standard Brownian motion and
	\begin{eqnarray}
		E_s&\equiv& \Big\{\exists s\in [q,t] y_1(tA(s/t))> 2\s tA(s/t)+(tA(s/t))^\g, \nonumber\\
		&&\qquad y_1(tA(q/t)) < 2\s t+B\sqrt{t}+R, y_1(t)<2\s t+B\sqrt{t}\Big\}
	\end{eqnarray}
	Let
	\begin{equation}\label{B3.27}
		\xi(s)= y_1(s)-\frac{s}{t-tA(q/t)}y_1(t) -\frac{t-s}{t-tA(q/t)} y_1(tA(q/t)),
	\end{equation}
	be a Brownian bridge from $0$ to $0$ in time $t-tA(q/t)$. Then the last expectation in \eqref{B3.26}, similarly as in \eqref{B3.24}, is bounded from above by
	\begin{multline}\label{B3.28}
		\E\left( \eee^{\l( y_1(t)-y_1(tA(q/t)))} \1_{\exists s\in [0,t-tA(q/t)]: \xi(s)> (s+tA(q/t))^\g -B\sqrt{t}-R }
		\right) \E\left( \eee^{\l (y_1(t)-y_1(tA(q/t)))}\right) 
		\\
		\times  \P\left(\exists s\in [0,t-tA(q/t)]: \xi(s)> (s+tA(q/t))^\g -B\sqrt{t}-R \right),
	\end{multline}
	where we used independence. By standard Brownian bridge estimates, the probability in \eqref{B3.28} converges to $0$ as $t$ tends to infinity. Overall, we can bound \eqref{B3.19} from above by
	\begin{eqnarray}
		\label{B3.29}
		\e/C K \eee^{-(1+2\s)t}\int_{t-r}^{t}\dd q \eee^{-\t^2(1-\rho^2)\left(t-tA(q/t)\right)}
		\eee^{2t-q} \E\left( \eee^{\l y_1(t-t A(q/t)) +\overline{\l} y_2(t-t A(q/t)) }\right)\nonumber
		\\
		\times \E\left(\eee^{2 \s x_1(tA(q/t))}\right) \leq \e,
	\end{eqnarray}
	for some large enough constant $C$ (which we can choose). To see this, we note that the left-hand side in \eqref{B3.29} coincides with the second moment (see \eqref{B3.2}) up to the $\e/C$ pre-factor. The second moment is finite by Lemma \ref{LemB3.1}. This proves Lemma \ref{B3.3}.
\end{proof}
Define
\begin{align}
	N_{\t,\s}^{B,\g}(t):= \sum_{l=1}^{n(t)}&e^{-t(1/2+\s^2)}\eee^{\s x_l(t)+i\t y_l(t)}\mathbbm{1}_{x_l(t)<2\s t+B\sqrt{t}, x_l\in \mathcal{S}_{r,t,\g}},
\end{align}
and let
\begin{equation}\label{eq:b_l}
	b_{l}(r):=\eee^{-r/2 -tA(r/t)\s^2}\eee^{\s x_l(r)+i\t x_l(r)}.
\end{equation}
The following lemma provides moment asymptotics for $N_{\t,\s}^{B,\g}(t)$ as $t\rightarrow \infty$.
\begin{lemma}[Moment asymptotics]\TH(Lem.const3) 
	For $\b\in B_3$, for any $B>0$ and $\g>1/2$,
	\be\Eq(con.100)
	\lim_{t\to\infty}\E\left[\left\vert N_{\t,\s}^{B,\g}(t)\right\vert^2\right]=C_{2,B}
	,
	\ee
	with $\lim_{B\to\infty} C_{2,B}=C_2$ and, for $k\in \N$, we have
	\be \Eq(con.101)
	\lim_{r\uparrow \infty}\lim_{t\to \infty} \E\left[\left\vert N_{\t,\s}^{B,\g}(t)\right\vert^{2k} ~\vert~ \mathcal{F}_r\right] =k! (C_{2,B}\MM_{2\s \s_b,0})^k \quad \mbox{a.s. and in }L^1.
	\ee
	Moreover, for $k'<k$,
	\be\Eq(con.102)
	\lim_{r\uparrow \infty}\lim_{t\to \infty} \E\left[N_{\t,\s}^{B,\g}(t)^k\overline{ N_{\t,\s}^{B,\g}(t)}^{k'} ~\big\vert~ \mathcal{F}_r \right] =0\quad \mbox{a.s. and in }L^1.
	\ee
	
\end{lemma}
For the proof of Lemma~\ref{Lem.const3} we need the following lemma.
\begin{lemma}\cite[Lemma~3.5]{hartklim18}\TH(Lem.help)
	Let $x,y$ be $\mathcal{N}(0,q)$ distributed random variables. Then, for any $m_1,m_2\geq 1$, and any constant $C>0$,
	\begin{align}\Eq(ny.11)
		\lim_{q \to \infty}\E&\left[ \left(\eee^{(m_1+2)\sigma x + i \tau m_2 x} +\eee^{(m_1+2)\sigma x - i\tau m_2 x}\right)  \1_{A_{x,\s,q,\gamma}}\right]\nonumber
		\\
		& o\Big(\eee^{2\s q} \E\left[ \left(\eee^{m_1\sigma x +i\tau m_2 x} +\eee^{m_1\sigma x -i\tau m_2 x}\right) \1_{A_{x,\s,q,\gamma}}\right]\E\left[\eee^{2\sigma y}  \1_{A_{y,\s,q,\gamma}}\right]\Big),
	\end{align}
	where $A_{x,\s,q,\gamma}:= \{x < 2\s q + Cq^{\gamma}\}.$
	Similarly,
	\begin{align}\Eq(ny.30)
		&\lim_{q \to \infty}\E\left[ \left(\eee^{(m_1+1)\sigma x +i\tau (m_2+1) x}+\eee^{(m_1+1)\sigma x -i\tau (m_2+1) x} \right)  \1_{A_{x,\s,q,\gamma}}\right]
		\\
		&= o\Big(\eee^{2\s q} \E\left[ \left(\eee^{m_1\sigma x +i\tau m_2 x} + \eee^{m_1\sigma x -i\tau m_2 x} \right)\1_{A_{x,\s,q,\gamma}}\right]\E\left[ \left(\eee^{(\sigma+i\tau) y} +\eee^{(\sigma-i\tau) y} \right) \1_{A_{y,\s,q,\gamma}}\right]\Big).\nonumber
	\end{align}
\end{lemma}
\begin{proof}[Proof of Lemma~\ref{Lem.const3}:]
	We proceed by induction over $k \in \N$.
	For $k=1$, we observe that, for $l\leq n(t)$,
	\begin{multline}
		\Eq(sun.1)
		1
		=
		\1\{x_l(t)>2\s t+B\sqrt{t}\} 
		+ \1\{x_l(t)<2\s t+B\sqrt{t}, x_l\not\in \mathcal{S}_{r,t,\g}\}
		\\
		+ \1\{x_l(t)<2\s t+B\sqrt{t}, x_l\in \mathcal{S}_{r,t,\g}\}
		.
	\end{multline}
	Plugging this decomposition of unity into \eqv(B3.2), we can rewrite $\E\left[\left\vert N_{\t,\s}(t)\right\vert^2\right]$ as
	\be\Eq(sun.2)
	\begin{aligned}
		\eee^{-2t\left(1/2+\s^2\right)}\E\Big[ 
		&\sum_{l_1,l_2=1}^{n(t)} \eee^{\s \left(x_{l_1}(t)+x_{l_2}(t)\right)+i\t \left(y_{l_1}(t)-y_{l_2}(t)\right)}\Big(\1_{\{x_{l_1}(t),x_{l_2}(t)>2\s t+B\sqrt{t}\}}
		\\
		&+ 2\1\{x_{l_1}(t)>2\s t+B\sqrt{t}, x_{l_2}(t)<2\s t+B\sqrt{t}\}
		\\
		&+\1\{x_{l_1}(t),x_{l_2}(t)<2\s t+B\sqrt{t}, x_{l_1}\not\in \mathcal{S}_{r,t,\g}\}
		\\
		&+ \1\{x_{l_1}(t),x_{l_2}(t) <2\s t+B\sqrt{t}, x_{l_2}\not\in  \mathcal{S}_{r,t,\g}
		\Big)\Big]
		\\
		&+\E\left[\left\vert N_{\t,\s}^{B,\g}(t)\right\vert^2\right]
		\\
		&=: \mathrm{(I)}+\mathrm{(II)}+\mathrm{(III)}+\mathrm{(IV)}+\mathrm{(V)}.
	\end{aligned}
	\ee
	Note that Terms~(I) and (III) can be made arbitrarily small by increasing $B$, resp. $r$ by computations as in the proofs of Lemmas~\thv(LemB3.2) and \thv(LemB3.3), respectively.
	
	Term~(II) can be treated as (I) as for the bounds in \eqv(B3.10) and below it suffices that one of the two particle positions exceeds $2 \s t +B\sqrt{t}$. 
	
	To control Term (IV), we upper bound it by $\eee^{-2t\left(1/2+\s^2\right)}$ times
	\begin{equation}\Eq(sun.3)
		\E\Big[ 
		\sum_{l_1,l_2=1}^{n(t)} \eee^{\s \left(x_{l_1}(t)+x_{l_2}(t)\right)+i\t \left(y_{l_1}(t)-y_{l_2}(t)\right)} 2\1\{x_{l_2}(t),x_{l_1}(t) <2\s t+B\sqrt{t}, x_{l_2}\not\in \mathcal{S}_{r,t,\g}\}\Big]
		.
	\end{equation}
	Observe that \eqv(sun.3) coincides with \eqv(B3.16). Following the argument
	after \eqv(B3.16), we see that \eqv(sun.3) can be made arbitrarily small by
	increasing $r.$
	
	Combining the bounds on the Terms (I), (II), (III) and (IV), the claim follows
	from Lemma~\thv(LemB3.1).
	
	To bound the $2k$-th moment, we rewrite \eqv(con.101) as
	\be\Eq(ny.1) 
	\begin{aligned}
		\frac{1}{2}\E\Big[\sum_{l_1,\dots, l_{2k}\leq n(t)}&\left(\prod_{j=1}^{2k}  \eee^{-t(1/2+\s^2 )} \eee^{\s x_{l_j}(t)+i\t y_{l_j}(t)}+ \prod_{j=1}^{2k}  \eee^{-t(1/2+\s^2 )} \eee^{\s x_{l_j}(t)-i\t y_{l_j}(t)}\right) 
		\\
		&
		\times \1\{x_{l_j}(t)<2\s t+B\sqrt{t},  x_{l_j} x_{l_1}\in \mathcal{S}_{r,t,\g}\}\Big] 
		,
	\end{aligned}
	\ee
	by grouping each summand together with its complex conjugate. 
	For $l_1,\dots, l_{2k}\leq n(t)$, we can find a matching using the following algorithm: Set $m=1$.
	\begin{itemize}\label{Optimal}
		\item[1.] Choose the two labels $j,j'$ such that $d(x_{l_j},x_{l_{j'}})$ is maximal. Call them $l_m$ and $l_{\sigma(m)}$ from know on.
		\item[2.] Delete them. Increase $m$ by $1$. Iterate.
	\end{itemize}
	We refer to the above algorithm as ``\textit{optimal matching}''. The pairs obtained in this way we denote by $(l_1,l_{\s(1)}), \dots, (l_k, l_{\s(k)})$. 
	We rewrite \eqv(ny.1) as
	\be\Eq(ny.2)
	\begin{aligned}
		\frac{1}{2}\E\Bigg[&\sum_{l_2,\dots, l_{k}\leq n(t)}\Big(\prod_{j=2}^{k} \eee^{-t(1+2\s^2 )} \eee^{\s \left(x_{l_j}(t)+x_{l_{\s(j)}}(t)\right)+i\t \left(y_{l_j}(t)+y_{l_{\s(j)}}(t)\right)}
		\\
		&   \qquad\qquad\qquad\quad \times\eee^{-t(1+2\s^2 )} \eee^{\s (x_{l_1}(t)+x_{l_{\s(1)}}(t))+i\t (y_{l_1}(t)-y_{l_{\s(1)}}(t))}
		+
		\\
		& \qquad\qquad\qquad
		\prod_{j=2}^{k} \eee^{-t(1+2\s^2 )} \eee^{\s \left(x_{l_j}(t)+x_{l_{\s(j)}}(t)\right)-i\t \left(y_{l_j}(t)+y_{l_{\s(j)}}(t)\right)}
		\\
		&  \qquad\qquad\qquad\quad \times \eee^{-t(1+2\s^2 )} \eee^{\s (x_{l_1}(t)+x_{l_{\s(1)}}(t))-i\t (y_{l_1}(t)-y_{l_{\s(1)}}(t))}
		\Big)
		\\
		& \quad \times
		\1\left\{x_{l_{\s(j)}}(t),x_{l_j}(t)<2\s t+B\sqrt{t},  x_{l_{\s(j)}}, x_{l_j}(s)\in  \mathcal{S}_{r,t,\g}\right\} 
		\\
		& \quad \times \1\left\{x_{l_{\s(1)}}(t),x_{l_1}(t)<2\s t+B\sqrt{t}, x_{l_{\s(1)}},x_{l_1}\in \mathcal{S}_{r,t,\g}\right\}
		\Bigg]
		.
	\end{aligned}
	\ee
	Using \eqv(cor.1), we can rewrite for $j\in\{1,\s(1)\}$    
	\begin{align}
		\Eq(ny.3)
		y_{l_j}(t)= \rho y_{l_j}(t)+ \sqrt{1-\rho^2} z_{l_j}(t),
	\end{align}
	where $(z_k(t))_{k\leq n(t)}$ are particles of a BBM on the same Galton-Watson tree as $(x_k(t))_{k\leq n(t)}$ but independent from it. Observe that using  the requirement that $d(x_{l_1},x_{l_{\s_1}})$ is chosen maximal, we have  
	\begin{multline}\Eq(ny.4)
		i\t (y_{l_1}(t)-y_{l_{\s(1)}}(t))=i\sqrt{1-\rho^2} \tau \left[z_{1}(t-tA(d(x_{l_1}(t),x_{l_{\s_1}}(t))/t))\right.\\
		\left. -z_{2}(t-tA(d(x_{l_1}(t),x_{l_{\s_1}}(t))/t))\right]
		+i \tau \rho \left(x_{l_1}(t)-x_{l_{\s(1)}}(t)\right),
	\end{multline}
	where $z_1, z_2$ are two independent $\mathcal{N}(0,t-tA(d(x_{l_1}(t),x_{l_{\s_1}}(t))/t)) $-distributed random variables.
	Plugging \eqv(ny.4) into \eqv(ny.2) and computing the expectation with respect to $z_1,z_2$ of the first summand (noting that the second is just its complex conjugate), we obtain
	\be\Eq(ny.5)
	\begin{aligned}
		&\frac{1}{2}\E\Bigg[\sum_{l_2,\dots, l_{k}\leq n(t)}\prod_{j=2}^{k} \eee^{-t(1+2\s^2)} \exp\left(\s (x_{l_j}(t)+x_{l_{\s(j)}}(t))+i\t (y_{l_j}(t)+y_{l_{\s(j)}}(t))\right)
		\\
		& \quad\quad \times \1\left\{x_{l_{\s(j)}}(t),x_{l_j}(t)<2\s t+B\sqrt{t}, x_{l_{\s(j)}}, x_{l_j}(s)\in  \mathcal{S}_{r,t,\g}\right\} 
		\\
		&
		\quad\quad \times \eee^{-t(1+2\s^2 )-\t^2(1-\rho^2)\left(t-tA(d(x_{l_1}(t),x_{l_{\s_1}}(t))/t)\right)} \eee^{(\s +i\tau\rho)x_{l_1}(t)+(\s -i\tau\rho)x_{l_{\s(1)}}}
		\\
		&
		\quad\quad\times \1\left\{x_{l_{\s(1)}}(t),x_{l_1}(t)<2\s t+B\sqrt{t}, x_{l_{\s(1)}}, x_{l_1}(s)\in  \mathcal{S}_{r,t,\g}\right\}
		\Bigg]
		.
	\end{aligned} 
	\ee
	We decompose 
	\be\Eq(ny.6)
	\begin{aligned}
		x_{l_{\s(1)}}(t)= x_{l_1}(d(x_{l_1},x_{l_{\s(1)}})) + x^{(1)}(t-tA(d(x_{l_1},x_{l_{\s(1)}}))/t); 
		\\
		x_{l_1}(t)= x_{l_1}(d(x_{l_1},x_{l_{\s(1)}})) + x^{(2)}(t-tA(d(x_{l_1},x_{l_{\s(1)}})/t)),
	\end{aligned}
	\ee
	where $x^{(1)},x^{(2)}$ are two independent $\mathcal{N}(0,t-tA(d(x_{l_1},x_{l_{\s(1)}})/t)) $-distributed random variables. By Step One of our matching procedure, we can plug \eqv(ny.5) into \eqv(ny.6) and compute the expectation with respect to $x^{(1)}$ and $x^{(2)}$, we obtain that \eqv(ny.5) is bounded from above by\footnote{A corresponding lower bound also holds due to the second moment computation in Lemma~\thv(Lem.const3).} 
	\be \Eq(ny.7)
	\begin{aligned}
		&\frac{1}{2}\E\Big[\sum_{l_2,\dots, l_{k}\leq n(t)}\prod_{j=2}^{k} \eee^{-t(1+2\s^2 )} \eee^{\s \left(x_{l_j}(t)+x_{l_{\s(j)}}(t)\right)+i\t \left(y_{l_j}(t)+y_{l_{\s(j)}}(t)\right)}
		\\
		& \quad \times \1\left\{x_{l_{\s(1)}}(t),x_{l_1}(t)<2\s t+B\sqrt{t}, x_{l_{1}}, x_{l_1}(s)\in  \mathcal{S}_{r,t,\g}\right\}
		\\
		&  \quad \times \eee^{-t(1+2\s^2 )-\t^2\left(t-tA(d(x_{l_1},x_{l_{\s(1)}})/t)\right)+\s^2 \left(t-tA(d(x_{l_1},x_{l_{\s(1)}})/t)\right)} \eee^{2\s x_{l_1}(d(x_{l_1},x_{l_{\s(1)}}))}
		\\
		& \quad \times
		\1\{  \forall s \in [r,d(x_{l_1},x_{l_{\s(1)}})] \colon x_{l_{\s(1)}}(s),x_{l_1}(s)\leq2\s tA(s/t)+(tA(s/t))^\g\}
		\Big]
		.
	\end{aligned} 
	\ee
	We now introduce the event 
	\begin{equation}
		\Eq(ny.8)
		\mathcal{A}_r= \Big\{\exists s \in [r, d(x_{l_1},x_{l_{\s(1)}})], \exists j\in\{2,\dots,k,\s(2),\dots, \s(k)\} \colon 
		d(x_{l_1}, x_{l_{j}}) =s\Big\} 
		.
	\end{equation}
	We can rewrite  \eqv(ny.7) as
	\be \Eq(ny.10)
	\E\left[\ldots \times \1_{\mathcal{A}_r}\right] + \E\left[\ldots \times \1_{\mathcal{A}_r^c}\right] =: J_{\mathcal{A}_r}+    J_{\mathcal{A}_r^c}.
	\ee
	We  prove that the first summand is of a smaller order than the second one using Lemma \ref{Lem.help}. Consider $ J_{\mathcal{A}_r}$. 
	Consider the skeleton generated by the leaves $l_1,l_{\s(1)}, \dots, l_k, l_{\s(k)}$ of the Galton-Watson tree. By $\mbox{path}(\cdot)$ we denote the unique path (= sequence of edges) leading from the given leaf ``$\cdot$'' to the root of the tree. To each edge in the Galton-Watson tree, we associate the following number
	\be\Eq(ny.17)
	m(e) := \sum_{j\in\{1,\s(1), \dots, k, \s(k)\}} \1_{e  \subset \mbox{path}(l_j)}.
	\ee
	
	For $k,j \leq n(t)$, define 
	\begin{align}\label{eq:length}
		\mathrm{length}^{\mathrm{eff}}(x_k(t), x_j(t)) & := tA\left(\frac{ d(x_1(t), x_k(t))}{t}\right)-tA\left(\frac{d(x_1(t), x_j(t))}{t} \right), \quad t \in \R_+
		,
	\end{align}
	the effective length.
	\begin{lemma}
		\label{lem:l_j_star}
		Consider the path of $x_{l_1}(t)$.
		There exists $l_{j*}$ which satisfies the following conditions
		\begin{itemize}
			\item[\textup{(i)}]
			$m$ is constant between $d(x_1(t), x_{j^*-1}(t))$ and $d(x_1(t), x_{j^*}(t))$ and, moreover,
			\begin{align}
				\mathrm{length}^{\mathrm{eff}}(x_{l_{j^*-1}},x_{l_{j^*}})>2r.
			\end{align}
			
			\item[\textup{(ii)}] 
			$\sum_{i=1}^{{j*-1}} \mathrm{length}^{\mathrm{eff}}(x_{l_{i-1}},x_{l_{i}})< (\mathrm{length}^{\mathrm{eff}} (x_{l_{j^{*}-1}},x_{l_{j^{*}}}))^\gamma$,
			where $\mathrm{length}^{\mathrm{eff}}$ is defined in \eqref{eq:length}.
		\end{itemize}
	\end{lemma}
	\begin{proof}
		Such a $l_{j^{*}}$ exists for all $t>t_0(r)$ because there are
		at most $2k-2$ points, where $m$ it is allowed to change. Hence, there must be a time interval of length $>2r$ (for $t$ large enough) during which $m$ does not change its value. Observe that if 
		\begin{equation}
			\sum_{i=1}^{{j*-1}} \mathrm{length}^{\mathrm{eff}}(x_{l_{i-1}},x_{l_{i}})>(2r)^\gamma ,
		\end{equation} 
		then only Condition~(ii) on $\mathrm{length}^{\mathrm{eff}}(x_{l_{j^*-1}},x_{l_{j^*}})$ needs to be checked. Assume that  $l_1,\dotsc,l_j$ all do not satisfy (ii). Then,
		\be\Eq(lenght.100)
		\sum_{i=1}^{{j-1}} \mathrm{length}^{\mathrm{eff}}(x_{l_{i-1}},x_{l_{i}})\leq Cr^{\left(\frac{1}{\gamma}\right)j}
		.
		\ee
		As $i<2k-2$ and the total time is equal to $t$, there must exist $j$ such that 
		\begin{align}
			\mathrm{length}^{\mathrm{eff}}(x_{l_{j-1}},x_{l_{j}})>Cr^{\left(\frac{1}{\gamma}\right)j}, \quad \text{for } t>t_0(r),
		\end{align}  
		where $t_0(r)$ is sufficiently large.
	\end{proof}
	

	We call the value of $m$ on the path of $x_{l_1}(t)$ between $d(x_1(t), x_{j^*-1}(t))$ and $d(x_1(t), x_{j^*}(t))$ $m^*$. Let us use the shortcut $R=d(x_1(t), x_{j^*-1}(t))$ and let 
	\begin{align}
		\ell = \ell(j^*, t) :=  d(x_1(t), x_{j^*}(t))-d(x_1(t), x_{j^*-1}(t)).
	\end{align}
	Then, on the time interval $(R,R+l)$, $m$ takes the value $m^*$.
	Moreover, at time $R$ the minimal
	particle is a.s.~$\min_{k\leq n(R)}x_k(R)+C\sqrt{2RtA(R/t)}$ to $\infty$ almost surely, for some $C>0$.
	To see this, observe that by symmetry, Markov's inequality and a Gaussian tail bound, 
	\be
	\P(\min_{k\leq n(R)}x_k(R) <-C\sqrt{2RtA(R/t)})\leq e^{R-(C')^2 tA(R/t)}.
	\ee
	As $A$ satisfies Assumption \ref{weak} this is summable in $t$ (and also $R$ for $C$ large enough). 
	The claim now follows by Borel--Cantelli. Hence, we may work on the event $x_{l_{j^*}}(R)>-C\sqrt{2RtA(R/t)}$. Now,
	\be \Eq(ny.19)
	x_{l_{j^*}}(R+\ell)- x_{l_{j^*}}(R)<x_{l_{j^*}}(R+\ell) +C\sqrt{2tRA(R/t)}
	.
	\ee
	Since we compute an expectation conditional on  $x_{l_{j*}}(R+\ell) < 2 \s tA((R+\ell)/t)+ (tA((R+\ell)/t))^\gamma$,
	we obtain on this event
	\be \Eq(ny.20)
	x_{l_{j^*}}(R+\ell)- x_{l_{j^*}}(R)< 2 \s tA((R+\ell)/t)+ (tA((R+\ell)/t))^\gamma+C\sqrt{2tRA(R/t)}.
	\ee
	Due to our choice of $j^*$, we have $2 \s tA(R/t)+C\sqrt{2tRA(R/t)}< C' [tA((R+\ell)/t)-tA(R/t)]^\g$ for some positive constant $C'$.
	By taking the expectation with respect to $x_{l_{j^*}}(R+\ell)- x_{l_{j^*}}(R)$ only, we can extract from $J_{\mathcal{A}}+\overline{J}_{\mathcal{A}}$ the factor
	\begin{align}
		\Eq(ny.21)
		\begin{split}
			\E\Big[ \left(\eee^{(m^*\s+ i\tau m' )x_{l_{j^*}}(R+\ell)- x_{l_{j^*}}(R)}+\eee^{(m^*\s- i\tau m' )x_{l_{j^*}}(R+\ell)- x_{l_{j*}}(R)}\right)
			\1_{ E_{j^*,R}}\Big],
		\end{split}
	\end{align}
	where $E_{j^*,R}$ is equal to
	\bea
	x_{l_{j*}}(R+\ell)- x_{l_{j*}}(R)&<& 2\s\left(tA\left(\frac{R+l}{t}\right)-tA\left(\frac{R}{t}\right)\right)\nonumber\\ &&+ (C'+1) \left(tA\left(\frac{R+l}{t}\right)-tA\left(\frac{R}{t}\right)\right)^\gamma.
	\eea
	By Lemma~\thv(Lem.help), \eqref{ny.21} is  
	\be \Eq(ny.22)
	\begin{aligned}
		o\Big(\eee^{2\s \ell} & \E\left[\eee^{((m^*-2)\s+ i\tau m' )x_{l_{j^*}}(R+\ell)- x_{l_{j^*}}(R)}\1_{ E_{j^*,R}}\right]
		\E\left[ \eee^{2\s (x_{l_{j*}}(R+\ell)- x_{l_{j*}}(R))}\1_{E_{j^*,R}}\right]\Big),
	\end{aligned}
	\ee
	for $l$ large (which by assumption~(i) of the lemma on $l$ corresponds to $r$ large). Note that the quantity, inside the brackets in \eqv(ny.22), corresponds to the same expectation but where in the underlying tree $l_1, l_{\s_1}$ branched off before time $R$.
	
	Iteratively, that leads to 
	\be\Eq(ny.23)
	J_{\mathcal{A}_r} +\overline{J}_{\mathcal{A}_r}\underset{t,r \to \infty}{=} o( J_{\mathcal{A}_r^c}+\overline{J}_{\mathcal{A}^c_r})
	.
	\ee
	Since $k$ was chosen arbitrarily, we know that the main contribution to the $2k$-th moment comes from the term where $l_1,\dots, l_k $ have split before time $r$ for $r$ large enough.
	We condition on $\mathcal{F}_r$ and compute:
	\be\Eq(ny.24)
	\begin{aligned}
		&\frac{1}{2}\E\Big[\sum_{l_1,l_2,\dots, l_{k}\leq n(t)}\Big(\prod_{j=2}^{k} \eee^{-t(1+2\s^2 )} \eee^{\s (x_{l_j}(t)+x_{l_{\s(j)}}(t))+i\t (y_{l_j}(t)-y_{l_{\s(j)}}(t))}\\
		&\qquad\qquad\qquad\qquad\quad+\prod_{j=2}^{k} \eee^{-t(1+2\s^2 )} \eee^{\s (x_{l_j}(t)+x_{l_{\s(j)}}(t))-i\t (y_{l_j}(t)-y_{l_{\s(j)}}(t))}\Big)
		\\
		& \times \1\left\{ x_{l_{\s(j)}}(t),x_{l_j}(t)<2\s t+B\sqrt{t}, x_{l_{\s(j)}}, x_{l_j}(s)\in  \mathcal{S}_{r,t,\g}\right\} 
		\1\{\sup_{j,j'\leq k}d(l_j,l_{j'})<r\} ~\Big|~ \mathcal{F}_r\Big]
		\\
		& =
		\E\Big[\sum_{l_1,l_2,\dots, l_{k}\leq n(t)}\prod_{j=2}^{k}   b_{l_j}(r)\overline{b}_{l_{\s(j)} }(r) \E\left[\left(\left(N_{\t,\s}^{B,\gamma}(t-r)\right)^{(j)}\right)^{2}\right] ~\Big|~ \mathcal{F}_r\Big],
	\end{aligned}
	\ee
	where $b_{l_j}(r)$ is defined in \eqv(eq:b_l) and $\left(N_{\t,\s}^{B,\gamma}(t-r)\right)^{(j)}$ are i.i.d.\ copies of $N_{\t,\s}^{B,\gamma}(t-r)$.
	By our second moment computations (Case $k=1$), as mentioned at the beginning of this proof,
	\be \Eq(ny.25)
	\lim_{t\to \infty} \E\left[\left(\left((N_{\t,\s}^{B,\gamma}(t-r)\right)^{(j)}\right)^{2}\right] = C_{2,B}.
	\ee
	Moreover, by invariance under permutation (in the labelling procedure),
	\be\Eq(ny.26)
	\sum_{l_1,l_2,\dots, l_{k}\leq n(t)}\prod_{j=2}^{k}   b_{l_j}(r)\overline{b}_{l_{\s(j)}}(r) = k!\Big(\sum_{k=1}^{n(r)} \eee^{2\s x_{k}(r)-(r+2\s^2 tA(r/t))}\Big)^{k}
	.
	\ee
	Observe that $\sum_{k=1}^{n(r)} \eee^{2\s x_{k}(r)-(r+2\s^2 tA(r/t))}$ converges in probability to $\MM_{2\s \s_b,0}$ by Lemma \ref{Lem.B1.2}. This proves \eqv(con.101).
	
	The case $k'<k$ follows similarly. Take an optimal matching (according to the procedure described below \eqref{ny.1}) of the first $k'$
	particles. The other particles will not be matched. Take one $l_1$ that has not
	been matched. Along its path, we can again find the first macroscopic piece on
	which $m(\cdot)$ is constant. Applying Lemma~\ref{Lem.help}, we get that the
	contribution is the largest if
	$\max_{j\in{1,\dots,k',1,\dots,k}}d(l_1,l_{j})<R$, for $R$ large enough. Observe
	that
	\be\Eq(ny.31)
	\E\Big[\sum_{k=1}^{n(t)}\eee^{\s x_k(t) +i\tau z_k(t)-(\frac{1}{2}+\s^2)t} ~\Big|~ \mathcal{F}_{R}\Big]=\eee^{(1-\s^2-\t^2)t/2-R-(\s^2-\t^2)tA(R/t))/2}\sum_{k=1}^{n(R)}\eee^{\s x_k(R) +i\tau z_k(R)}.
	\ee
	We note that by Taylor expansion $tA(R/t)=\s_bR+O(R^2/t)$. As  $1-\s^2-\t^2<0$ in $B_3$, the sum on the r.h.s.~of \eqv(ny.31) converges to zero as $t\uparrow\infty$.
	This together with the argument in the even case implies Lemma~\thv(Lem.const3).
\end{proof}

\begin{proof}[Proof of Theorem~\thv(clt.B3)]
	
	Recall that the even (resp., odd) moments of the complex isotropic distribution $\mathcal{N}(0,C_{2,B}\MM_{2\s\s_b,0})$ coincide with the r.h.s.\ of
	\eqref{con.101} (resp., \eqref{con.102}).
	By Lemma~\thv(Lem.const3), conditionally on $\mathcal{F}_r$, the moments of
	$N_{\s,\t}^{c,B}(t)$ converge to the moments of a $\mathcal{N}(0,C_{2,B}\MM_{2\s,0})$ a.s.\ as $t\uparrow \infty$ and then $r\uparrow
	\infty$. Since the
	normal distribution is uniquely characterised by its moments, this implies
	convergence in distribution.  Moreover, by Lemma~\thv(LemB3.2) and
	Lemma~\thv(LemB3.3),
	\be \Eq(con.501)
	\lim_{A\uparrow\infty} \lim_{t\uparrow \infty} \mathcal{L}\left[N_{\s,\t}(t)-N_{\s,\t}^{c,B}(t) \right] = \delta_0, 
	\ee
	and $\lim_{B\to\infty} C_{2,B}=C_2$. The claim of Theorem~\thv(clt.B3) follows.
\end{proof}

\section{Proof of Theorem~\ref{Cor:phase-diagram}}
\label{sec:phase-diagram}

In this section, as a consequence of the fluctuation results of the previous sections, we
derive the phase diagram shown on Fig.~\ref{fig-rem-phase-diagram}.
\begin{proof}[Proof of Theorem~\ref{Cor:phase-diagram}]
	Convergence in probability for $\beta \in  B_1$ and $ B_3$   in
	\eqref{eq:limiting-log-partition-function} follows from Theorems~\ref{paseB1} and
	\ref{clt.B3} by \cite[Lemma~3.9~(1)]{KaKli14}. Convergence for the glassy
	phase $\beta \in \overline{B_2}$ follows from Theorem~\ref{phaseB2.1} for $|\rho|=1$ and from Theorem~\ref{phaseB2.2} otherwise. The formula \eqref{eq:limiting-log-partition-function} for the boundaries of the phases follows from the continuity of the limiting log-partition function.
\end{proof}

\bibliographystyle{abbrv}
\bibliography{complexcrem-ref}

\end{document}